\documentclass[11pt]{amsart}
\usepackage[utf8]{inputenc}
\usepackage[english]{babel}

\usepackage{amscd,amssymb}
\usepackage{amsmath}
\usepackage{amsthm}
\usepackage{graphicx}
\usepackage{fancybox}
\usepackage{epic,eepic}
\usepackage{amstext}
\usepackage{mathtools}
\usepackage{esint}
\usepackage[square,numbers]{natbib}
\usepackage{booktabs}
\usepackage[hide links]{hyperref}
\usepackage{comment}
\usepackage{mathtools}
\usepackage{tikz,pgfplots}
\usepackage{subcaption}

\def\V{\mathcal{V}}

\def\W{\mathcal{W}}

\def\Pnull{\mathbb{P}^0}

\def\elem{T}

\def\T{\mathcal{T}}
\def\TH{\mathcal{T}_H}
\def\l{\ell}

\def\C{\mathcal{C}}

\def\PiH {\Pi_H}

\def\PitH{\widetilde{\Pi}_H}
\def\diam{\mathrm{diam}}
\def\linh{\mathrm{span}}
\def\supp{\mathrm{supp}}

\def\one{\mathbf{1}}

\newcommand{\tnormf}[2]{{\| #1 \|}_{L^2(#2)}}

\newcommand{\vnorm}[2]{{\left\lVert #1 \right\rVert}_{a,#2}}
\newcommand{\vnormf}[2]{{\| #1 \|}_{a,#2}}

\newcommand{\tspf}[3]{{( #1\,,\,#2 )}_{L^2(#3)}}

\newcommand\dx{\,\text{d}x}

\DeclareMathOperator*{\argmin}{argmin}

\definecolor{myBlue}{RGB}{113,104,238} 
\definecolor{myGreen}{RGB}{114,175,30} 
\definecolor{myRed}{RGB}{180,50,50}  
\definecolor{myOrange}{RGB}{225,92,22}

\def\dx{\,\text{d}x}

\newtheorem{theorem}{Theorem}[section]

\newtheorem{lemma}[theorem]{Lemma}
\newtheorem{assumption}[theorem]{Assumption}

\theoremstyle{definition}

\theoremstyle{remark}
\newtheorem{remark}[theorem]{Remark}
\numberwithin{theorem}{section}
\numberwithin{equation}{section}
\numberwithin{table}{section}
\numberwithin{figure}{section}

\allowdisplaybreaks

\begin{document}
	
\title[]{Super-localization of elliptic multiscale problems}
\author[]{Moritz Hauck$^\dagger$, Daniel Peterseim$^{\dagger,\ddagger}$}
\address{${}^{\dagger}$ Institute of Mathematics, University of Augsburg, Universit\"atsstr.~12a, 86159 Augsburg, Germany}
\address{${}^{\ddagger}$ Centre for Advanced Analytics and Predictive Sciences (CAAPS), University of Augsburg
	Universit\"atsstr.~12a, 86159 Augsburg, Germany}
\email{\{moritz.hauck, daniel.peterseim\}@uni-a.de}

\thanks{The work of Moritz Hauck and Daniel Peterseim is part of a project that has received funding from the European Research Council (ERC) under the European Union's Horizon 2020 research and innovation programme (Grant agreement No.~865751 --  RandomMultiScales).}

\thanks{This paper is dedicated to Professor Carsten Carstensen on the occasion of his 60th birthday}
\maketitle


\begin{abstract}
Numerical homogenization aims to efficiently and accurately approximate the solution space of an elliptic partial differential operator  with arbitrarily rough coefficients in a $d$-dimensional domain.
The application of the inverse operator to some standard finite element space defines an approximation space with uniform algebraic approximation rates with respect to the mesh size parameter $H$. This holds even for under-resolved rough coefficients. However, the true challenge of numerical homogenization is 
the localized computation of a localized basis for such an operator-dependent approximation space. This paper presents a novel localization technique 
that leads to a super-exponential decay of its basis relative to $H$. This suggests that basis functions with supports of width $\mathcal O(H|\log H|^{(d-1)/d})$ are sufficient to preserve the optimal algebraic rates of convergence in $H$ without pre-asymptotic effects. A sequence of numerical experiments illustrates the significance of the new localization technique when compared to the so far best localization to supports of width $\mathcal O(H|\log H|)$.
\end{abstract}
\indent
{\tiny {\bf AMS subject classifications.}  
	{\bf 65N12}, 
	{\bf 65N30} 
} 
\section{Introduction}
This paper considers the numerical solution of second order elliptic partial differential equations (PDE) with a strongly heterogeneous and highly varying (non-periodic) coefficient and an $L^2$-right-hand side. The heterogeneity and oscillations of the coefficient may appear on several non-separated scales.  For this general class of coefficients, classical finite element methods based on universal, problem-independent polynomial ansatz spaces can perform arbitrarily badly; see e.g. \cite{BabO00}. In order to cure this problem, problem-specific information needs to be incorporated in the ansatz space. This approach is often referred to as numerical homogenization. For an overview on numerical homogenization beyond periodicity and scale separation, we refer to the review paper \cite{Peterseim2021} and the recent textbooks \cite{OwhS19,MalP20}. 

A common technique for the construction of problem-adapted ansatz spaces in numerical homogenization is the application of the inverse operator to some standard finite element spaces on a coarse mesh which does not necessarily resolve the coefficient. This technique has been used explicitly or implicitly in many works; see e.g. \cite{OwZ11,GGS12,MaP14,OZB14,Owh17}. It is straightforward that the corresponding Galerkin method yields optimal rates of convergence on arbitrarily coarse meshes without any pre-asymptotic effects. However, this approach is not feasible without modification because, in general, global problems need to be solved to compute the basis functions of the ansatz space which, in turn, tend to have global support. Hence, the true challenge of numerical homogenization is to identify a local basis of the ansatz space or at least one that consists of rapidly decaying basis functions. This rapid decay makes it possible to approximate the global basis functions by localized counterparts which are solutions to problems on local patches of the coarse grid. The corresponding Galerkin approximation in the space spanned by the localized basis then yields a practical numerical homogenization method. 

In early works on this localization problem, algebraically decaying basis functions were constructed; see e.g. \cite{OwZ07,Berlyand2010FluxNA,OwZ11}. 
A near-optimal localization was later achieved by the Localized Orthogonal Decomposition (LOD) method \cite{MaP14,HeP13,KPY18,BrennerLOD} which constructs a fixed number of basis functions per mesh entity which decay exponentially in modulus relative to the coarse mesh size $H$. Hence, computable local basis functions on patches of diameter $\mathcal O(H|\log H|)$ are sufficient to preserve optimal algebraic rates of convergence in $H$ without any pre-asymptotic effects. Other localization techniques such as the AL basis \cite{GGS12} and generalized finite element methods \cite{BaL11,Chen21,Ma21} lead to qualitatively similar results. 
The latter approaches are based on local spectral computations and the partition of unity method. Its basis functions are local by construction. However, algebraic accuracy in $H$ typically requires to moderately enlarge the number of basis functions per mesh entity. Hence, the dimension of the approximation space grows by a logarithmic-in-$H$ multiplicative factor. The number of non-zero entries in the resulting stiffness matrices leads to a rough complexity estimate of the different methods. Compared to the tentative minimum number of non-zero entries which is proportional to the number of mesh entities, the multiplicative overhead is $\mathcal O(|\log H|^{d})$ for LOD-type methods and $\mathcal O(|\log H|^{d+1})$ for the other techniques according to the existing theoretical predictions. Even in practice, none of the known constructions seem to reduce the overhead below $\mathcal O(|\log H|^{d})$ without reducing the approximation rate.

This paper introduces a novel localization strategy in the context of the LOD that reduces the diameter of the supports of the basis functions to $\mathcal O(H|\log H|^{(d-1)/d})$, in practice without loss of accuracy. On the level of the sparsity of the stiffness matrix, this corresponds to the favorable overhead of $\mathcal O(|\log H|^{d-1})$. 
The novel localization approach identifies local finite element source terms that yield rapidly decaying responses under the solution operator of the PDE. Localized ansatz functions are then obtained by solving patch-local problems with zero boundary data for the previously determined source terms. 

While classical LOD theory \cite{MaP14,HeP13,Peterseim2021} only proves that the novel localization approach at least recovers classical LOD performance, an even super-exponentially decay of the localization error is conjectured. Although a rigorous proof is still open, we provide a justification employing a conjecture related to spectral geometry. Numerical experiments in two and three spatial dimensions confirm the super-exponential decay and demonstrate the computational superiority of the novel localization strategy.

The remaining part of the paper is outlined as follows. Section \ref{s:mp} recalls the prototypical elliptic model problem and some basics from finite element theory. A prototypical numerical homogenization method is presented in Section \ref{s:idealmethod}. Using a novel localization approach, this method is then turned into a feasible practical method in Sections \ref{s:locapprox} and \ref{s:locmethod}. Section~\ref{s:error} presents an abstract a-priori error analysis and Section~\ref{s:qo} justifies the super-exponential decay of the localization error. Numerical experiments are presented in  Section~\ref{s:ne} and the paper closes with a conclusion and an outlook on future research.

\section{Model Problem}\label{s:mp}
Let us consider the prototypical second order elliptic PDE $-\mathrm{div}A\nabla u = f$ in weak form with homogeneous Dirichlet boundary conditions on a polygonal Lipschitz domain $\Omega\subset \mathbb{R}^d$, $d\in\mathbb{N}$. Without loss of generality, we assume that $\Omega$ is scaled to unit size, i.e., its diameter is of order one. The matrix-valued coefficient function $A\in L^\infty(\Omega,\mathbb{R}^{d\times d})$ is supposed to be symmetric and positive definite almost everywhere. More precisely, we assume that there exist constants $0<\alpha\leq \beta< \infty$ such that for almost all $x \in \Omega$ and for all $\eta\in \mathbb{R}^d$
\begin{equation}\label{eq:propA}
\alpha|\eta|^2 \leq (A(x)\eta)\cdot \eta \leq \beta |\eta|^2
\end{equation}
with $|\cdot|$ denoting the Euclidean norm of a $d$-dimensional vector.
The solution space is the Sobolev space $\V \coloneqq H^1_0(\Omega)$ and the bilinear form  $a\colon\V\times\V\rightarrow \mathbb{R}$ associated with the problem is given by
\begin{equation*}
a(u, v) \coloneqq
\int_\Omega (A\nabla u)\cdot\nabla v\dx.
\end{equation*}
The symmetry and the condition \eqref{eq:propA} ensure that the above bilinear form is an inner product on $\V$. Its induced norm is the energy norm  $\vnorm{\cdot}{\Omega} \coloneqq \sqrt{a(\cdot , \cdot)}$ which is
equivalent to the canonical Sobolev norm on $\V$. Given $f\in L^2(\Omega)$, the unique weak solution $u \in \V$ of the boundary value problem satisfies, for all $v \in \V$, 
\begin{equation*}
a(u, v) = (f,v)_{L^2(\Omega)}.
\end{equation*}
Note that the moderate restriction to right-hand sides in $L^2(\Omega)$ (rather than the dual space  $\V^\prime=H^{-1}(\Omega)$) will be source of  the uniform linear convergence of the numerical homogenization method although the possibly rough coefficient in general prevents $H^2(\Omega)$-regularity which would be required by classical finite elements. We refer to \cite{Peterseim2021} for a detailed discussion and possible generalizations to right-hand sides with less regularity. Henceforth, we will refer to $\mathcal A^{-1}\colon L^2(\Omega)\rightarrow \V$ as the solution operator that maps $f\in L^2(\Omega)$ to the unique solution $u\in\V$ of the above weak formulation.

\section{Optimal operator-dependent approximation}\label{s:idealmethod}
This section introduces prototypical operator adapted ansatz spaces computed by applying the inverse operator to right-hand sides in classical finite element spaces. This approach transforms the approximation problem of the solution into an approximation problem of the right-hand side by classical finite element spaces.  For this purpose, we shall introduce the, possibly coarse, quasi-uniform mesh $\TH$ which is a finite subdivision of $\Omega$ into closed, convex, and shape regular elements with diameter at most $H>0$.
This characterization of admissible meshes includes both, structured quadrilateral or Cartesian meshes and unstructured shape-regular simplicial meshes as they are well-established in the theory of finite elements.

For the approximation of $f\in L^2(\Omega)$, we consider the simple choice of $\T_H$-piecewise constant functions 
$$\Pnull(\T_H)=\operatorname{span}\{\one_\elem\,|\,\elem\in \mathcal T_H\}$$ spanned by the basis of the characteristic functions $\one_\elem$ of the mesh elements $\elem\in\TH$. By using different finite element spaces, it is possible to construct, e.g., higher order methods \cite{Mai20ppt} or multi-level methods \cite{Owh17,FeP20,HaPe21}.
Let $\Pi_H\colon L^2(\Omega)\rightarrow \Pnull(\TH)$ denote the $L^2$-orthogonal projection onto $\Pnull(\TH)$ and recall that, for all $\elem\in \mathcal T_H$, it satisfies the following (local) stability and approximation properties
\begin{equation}\label{e:L2proj}
\begin{aligned}
\|\PiH  v\|_{L^2(\elem)}&\leq \|v\|_{L^2(\elem)}\qquad &&\text{ for all }v\in L^2(\elem),\\
\|v-\PiH v\|_{L^2(\elem)}&\leq {\pi}^{-1}H \|\nabla v\|_{L^2(\elem)}\quad &&\text{ for all }v\in H^1(\elem).
\end{aligned}
\end{equation}

Given $\Pnull(\T_H)$, the prototypical operator-adapted ansatz space $\V_H$ is obtained by applying the solution operator, i.e.,  
\begin{equation}\label{eq:ideallod}
\V_H\coloneqq\linh\{\mathcal{A}^{-1}\one_\elem\,|\,\elem\in \mathcal T_H\}.
\end{equation}
The Galerkin method based on the operator-adapted ansatz space $\V_H$ seeks a discrete approximation $u_H\in \V_H$ such that, for all $v_H\in \V_H$,
\begin{equation}\label{e:galerkinideal}
a(u_H,v_H) = \tspf{f}{v_H}{\Omega}.
\end{equation}
Let $\mathcal A_H^{-1}\colon  L^2(\Omega)\rightarrow \V_H$, $f\mapsto u_H$ denote the discrete solution operator of the above Galerkin method. 

\begin{remark}[Exactness for $\T_H$-piecewise constant right-hand sides]\label{r:exactness}
	The discrete solution operator $\mathcal A_H^{-1}$ is exact for $f \in \Pnull(\TH)$. This implies that $\mathcal A_H^{-1}\circ \PiH = \mathcal A^{-1}\circ \PiH$.
\end{remark}
The following lemma states an approximation result on the operator level without pre-asymptotic effects. 

\begin{lemma}[Uniform operator approximation]\label{l:ua}
	For any $s\in [0,1]$, the discrete solution operator $\mathcal A^{-1}_H$ uniformly approximates  $\mathcal A^{-1}$ in the operator norm $\|\cdot \|_{H^s(\Omega)\rightarrow \V}$ with an algebraic rate in $H$, i.e., there is $C>0$ independent of $\alpha,\beta,H$ such that, for all $f\in H^s(\Omega)$,
	\begin{equation*}
	\pi \alpha^{1/2}\,{\vnormf{\mathcal A^{-1}f - \mathcal A_H^{-1}f}{\Omega}} \leq H \tnormf{f -\PiH f}{\Omega}\leq  C H^{1+s}{\|f\|_{H^s(\Omega)}}.
	\end{equation*}
\end{lemma}
\begin{proof}
	This is a well-known result; see e.g. \cite{MalP20,Peterseim2021}. A short proof is presented for completeness. Define the error $ e\coloneqq\mathcal A^{-1}f -\mathcal A_H^{-1} f$. Since $a(u_H,e) = 0$ by symmetry and Galerkin orthogonality, we obtain
	\begin{align*}
	\vnormf{ e}{\Omega}^2&= a( e, e) = a(u,e) = \tspf{f-\PiH f}{ e}{\Omega}\leq \tnormf{f-\PiH f}{\Omega}\tnormf{ e}{\Omega}
	\end{align*}
	using that Remark \ref{r:exactness} implies for all $\elem \in \TH$
	$$\tspf{\one_\elem}{e}{\Omega} = a(\mathcal A^{-1}\one_\elem,e) = a(\mathcal A_H^{-1}\one_\elem,e) = 0,$$  i.e., $\PiH e = 0$. This and \eqref{e:L2proj} yield
		\begin{align*}
	\tnormf{ e}{\Omega}&=\tnormf{ e-\PiH e}{\Omega}\leq \pi^{-1}\alpha^{-1/2}  H  \vnormf{ e}{\Omega}.
	\end{align*}
The combination of the previous estimates and \eqref{e:L2proj} readily yields the result for $s\in\{0,1\}$. For the other values of $s$, the assertion can be concluded by arguments from interpolation theory (see e.g. \cite{BrS08}) introducing a constant $C>0$ independent of $A$ and $H$. 
\end{proof}

\section{Novel localization strategy}\label{s:locapprox}
The canonical basis functions  $\{\mathcal A^{-1}\one_\elem\,|\,\elem \in \TH\}$ of the operator-adapted approximation space $\V_H$, defined in \eqref{eq:ideallod}, are non-local and have a slow (algebraic) decay driven by the decay of the Green's function associated with the PDE operator $\mathcal A$. By a localization of the basis, it is possible to deduce a practically feasible variant of method \eqref{e:galerkinideal}. State-of-the-art localization strategies based on the decay of the fine-scale Green's function \cite{MaP14,HeP13} or on subspace decomposition \cite{KorY16,KPY18,BrennerLOD} yield an exponential decay of the localization error as layers of neighboring elements are added. The fast decay is triggered by a Lagrange property of the basis with respect to some quasi-interpolation operator on the coarse mesh which corresponds to $\Pi_H$ in the present setting. This section presents an advanced localization strategy which even promises super-exponentially decaying localization errors.

The novel localization approach identifies local $\TH$-piecewise constant source terms that yield rapidly decaying (or even local) responses under the solution operator $\mathcal A^{-1}$ of the PDE. As usual, localization relies on the concept of (local) patches in the coarse mesh $\T_H$ based on neighborhood relations between mesh elements. Given a union of elements $S\subset \Omega$,
the first order element patch $\mathsf{N}(S)=\mathsf{N}^1(S)$ of $S$ is given by
\begin{equation*}
\mathsf{N}^1(S)\coloneqq \bigcup \left\{\elem\in \TH\,|\,\elem\cap S\neq \emptyset\right\}.
\end{equation*}
For any $\ell=2,3,4,\dots$, the $\ell$-th order patch $\mathsf{N}^\ell(T)$ of $T$ is then given recursively by
\begin{equation*}
\mathsf{N}^\ell(T)\coloneqq \mathsf{N}^{1}(\mathsf{N}^{\ell-1}(T)).
\end{equation*}

In order to simplify the notation in the subsequent derivation, we shall fix an arbitrary element $\elem\in\T_H$ and the oversampling parameter $\ell\in\mathbb{N}$. We will refer to the $\ell$-th order patch of $T$ by $\omega\coloneqq \mathsf{N}^\ell(\elem)$  and make the meaningful assumption that $\omega$ does not coincide with the whole domain $\Omega$. Let  $\T_{H,\omega}$ denote the submesh of $\TH$ with elements in $\omega$ and let $\Pi_{H,\omega}\colon  L^2(\omega) \rightarrow \Pnull(\T_{H,\omega})$ denote  the $L^2(\omega)$-orthogonal projection onto $\Pnull(\T_{H,\omega})$.

The (ideal) basis function $\varphi = \varphi_{\elem,\ell}\in \V_H$ associated with the element $\elem$ is given by the ansatz
\begin{equation*}
\varphi = \mathcal A^{-1}g \quad \text{ with  }\quad g = g_{\elem,\ell} \coloneqq \sum_{K\in \mathcal T_{H,\omega}}c_K \one_K,
\end{equation*}
where $(c_K)_{K \in \T_{H,\omega}}$ are unknown coefficients to be determined subsequently. Given $\varphi$, its Galerkin projection onto the local subspace $H^1_0(\omega)$ defines a localized approximation $\varphi^\mathrm{loc} = \varphi^\mathrm{loc}_{\elem,\ell} \in H^1_0(\omega)$ that satisfies, for all $v\in H^1_0(\omega)$,
\begin{equation}\label{e:patchproblem}
a_\omega(\varphi^\mathrm{loc},v)\coloneqq\int_\omega (A\nabla \varphi^\mathrm{loc})\cdot\nabla v\dx  = \tspf{g}{v}{\omega}.
\end{equation}
Note that, throughout this paper, we do not distinguish between $H^1_0(\omega)$-functions and their  $\V$-conforming extension by zero to the full domain $\Omega$. 
In general, the local function $\varphi^\mathrm{loc}$ is a poor approximation of the possibly global function $\varphi$. However, the appropriate choice of $g$ will lead to a highly accurate approximation in the energy norm.

The choice of $g$ requires a quick reminder on traces of $H^1(\omega)$-functions (see e.g. \cite{LiM72a} for details). Let 
\begin{equation*}
\gamma_0 = \gamma_{0,\omega}\colon \V_ \omega \rightarrow X \coloneqq \mathrm{range}\,\gamma_{0}\subset H^{1/2}(\partial  \omega)
\end{equation*}
denote the trace operator on $\omega$ restricted to the complete subspace $\V_\omega \coloneqq  \{v\vert_\omega\,\vert\,v\in \V\}\subset H^1(\omega)$. The $A$-harmonic extension $\gamma_{0}^{-1}$ defines a continuous right-inverse of $\gamma_{0}$. Given $w \in X$,  $\gamma_{0}^{-1}w = \gamma_{0,\omega}^{-1}w$ satisfies $\gamma_0 \gamma_{0}^{-1}w= w$ and, for all $v \in H^1_0( \omega)$,
\begin{equation}
\label{eq:harmext}
a_\omega(\gamma_{0}^{-1}w,v)= 0.
\end{equation}

Given the trace and extension operators, we conclude from \eqref{e:patchproblem} and \eqref{eq:harmext}, using that $v - \gamma_0^{-1}\gamma_0 v \in H^1_0(\omega)$, that 
	\begin{equation*}
	a(\varphi^\mathrm{loc} ,v) = a_\omega(\varphi^\mathrm{loc} ,v) = a_\omega(\varphi^\mathrm{loc} ,v-\gamma_{0}^{-1} \gamma_{0}\, v) = (g,v- \gamma_{0}^{-1} \gamma_{0}\, v)_{L^2(\omega)}.
	\end{equation*}
	This, the definition of $\varphi$ and $\supp\, g\subset\omega$ yield the following crucial observation 
	\begin{equation}
	\label{e:nd1}
	a(\varphi-\varphi^\mathrm{loc} ,v) = (g,v)_{L^2(\omega)}-a_\omega(\varphi^\mathrm{loc} ,v)= (g, \gamma_{0}^{-1} \gamma_{0}\, v)_{L^2(\omega)}
	\end{equation}
	for any $v \in \V$. 

This observation rephrases the smallness of the localization error as the (almost) $L^2$-orthogonality of $g$ to the space
\begin{equation}\label{e:Y}
Y\coloneqq \gamma_{0}^{-1}X\subset \V_\omega
\end{equation}
of $A$-harmonic functions on $\omega$ (which satisfy the homogeneous Dirichlet boundary condition on $\partial \Omega\cap\partial\omega$).
	It shall be noted that \eqref{e:nd1} coincides with the definition of the conormal derivative $A\nabla \varphi^\mathrm{loc}\cdot n\in X^\prime$ ($n$ denoting the outer normal unit vector with respect to $\omega$), i.e., the smallness of the localization error is also equivalent to a small $X^\prime$-norm of the conormal derivative of $\varphi^\mathrm{loc}$.

An optimal choice of $g$ is, hence, realized by the singular value decomposition (SVD) of the operator $\Pi_{H,\omega}|_Y$	which has finite rank less or equal to $N \coloneqq \#\T_{H,\omega}$. The SVD ist given by 
	\begin{equation*}
	\Pi_{H,\omega}|_Y\,v = \sum_{k = 1}^{N} \sigma_k (v,w_k)_{H^1(\omega)}\,g_k
	\end{equation*}
	with singular values $\sigma_1\geq \dots \geq \sigma_{N}\geq 0$, $L^2(\omega)$-orthonormal left singular vectors $g_1\,\ldots,g_N$, and  $H^1(\omega)$-orthonormal right singular vectors $w_1\,\ldots,w_N$. The  left singular vector $g_N$ corresponding to the smallest singular value $\sigma_N$ satisfies
	\begin{equation}\label{e:gopt}
	g_N\in\argmin_{g\in \mathbb P^0(\T_{H,\omega})\colon\|g\|_{L^2(\omega)}=1}\;\sup_{v\in Y\colon \|v\|_{H^1(\omega)} = 1} {(g, v)_{L^2(\omega)}}.\end{equation}
	In this sense, $g_N$ is an optimal choice for $g$ and we define 
	\begin{equation}\label{e:orth}
\sigma_\elem = \sigma_\elem(H,\ell) \coloneqq \sigma_N =  \sup_{v\in Y\colon  \|v\|_{H^1(\omega)} = 1} {\tspf{g_N}{v}{\omega}}.
\end{equation}
The value $\sigma_\elem$ will be used in the remainder as a measure for the (quasi-) orthogonality between $g$ and $Y$. It coincides with the ${X^\prime}$-norm of the conormal derivative of $\varphi^\mathrm{loc}$ up to a constant depending only on the geometry of the patch.

\begin{figure}[h]
	\includegraphics[width=.49\linewidth]{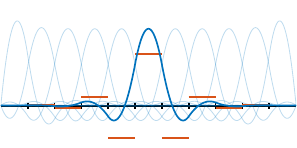}\hfill
	\includegraphics[width=.49\linewidth]{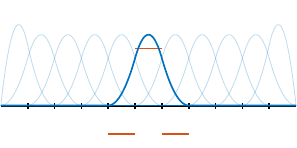}
	\caption{Classical (ideal) LOD basis (left) and novel local basis obtained by the novel localization approach (right) for $A \equiv 1$. The ($L^2$-normalized) right-hand sides $g$ corresponding to the central basis functions are depicted in orange.}
	\label{figslodlod}
\end{figure}

\begin{remark}[Local basis in 1d]\label{r:1d}
	For one spatial dimension, the trace space is at most two-dimensional and so is the space $Y$. Thus, already for $\ell = 1$, we can choose an $L^2$-normalized $g\in \Pnull(\T_{H,\omega})$ that is $L^2$-orthogonal to $Y$.  
	Hence, the novel localization approach yields a local basis.
	It is easy to verify that, for a constant coefficient $A$, the constructed basis coincides with the canonical quadratic B-splines; see also \cite[Chapter~2]{Bsplines}.
	In comparison, the ideal LOD basis functions \cite{MaP14,Peterseim2021} and its corresponding right-hand sides are globally supported but decay exponentially fast; see \cite[Theorem 3.15]{Peterseim2021} and Lemma \ref{l:L2flux}. By construction, the LOD basis functions fulfill a Lagrange-type property with respect to element-averages. For an illustration of the novel basis and the LOD basis, see Figure \ref{figslodlod}.
\end{remark}
The possible non-uniqueness of the smallest singular value becomes immediately  clear in the one-dimensional case, where, for any $\ell \geq 2$, a simple counting argument shows that there are multiple optimal choices of $g$. In higher dimensions, the problem rather manifests itself in clusters of small singular values which may appear for large $\ell$ and certain geometric configurations of the patches near the boundary. The choice of an appropriate $g$ can then be difficult; see Appendix \ref{s:impl} for an illustrative example and a practical solution.

\begin{figure}[h]
	\includegraphics[width=.3\linewidth]{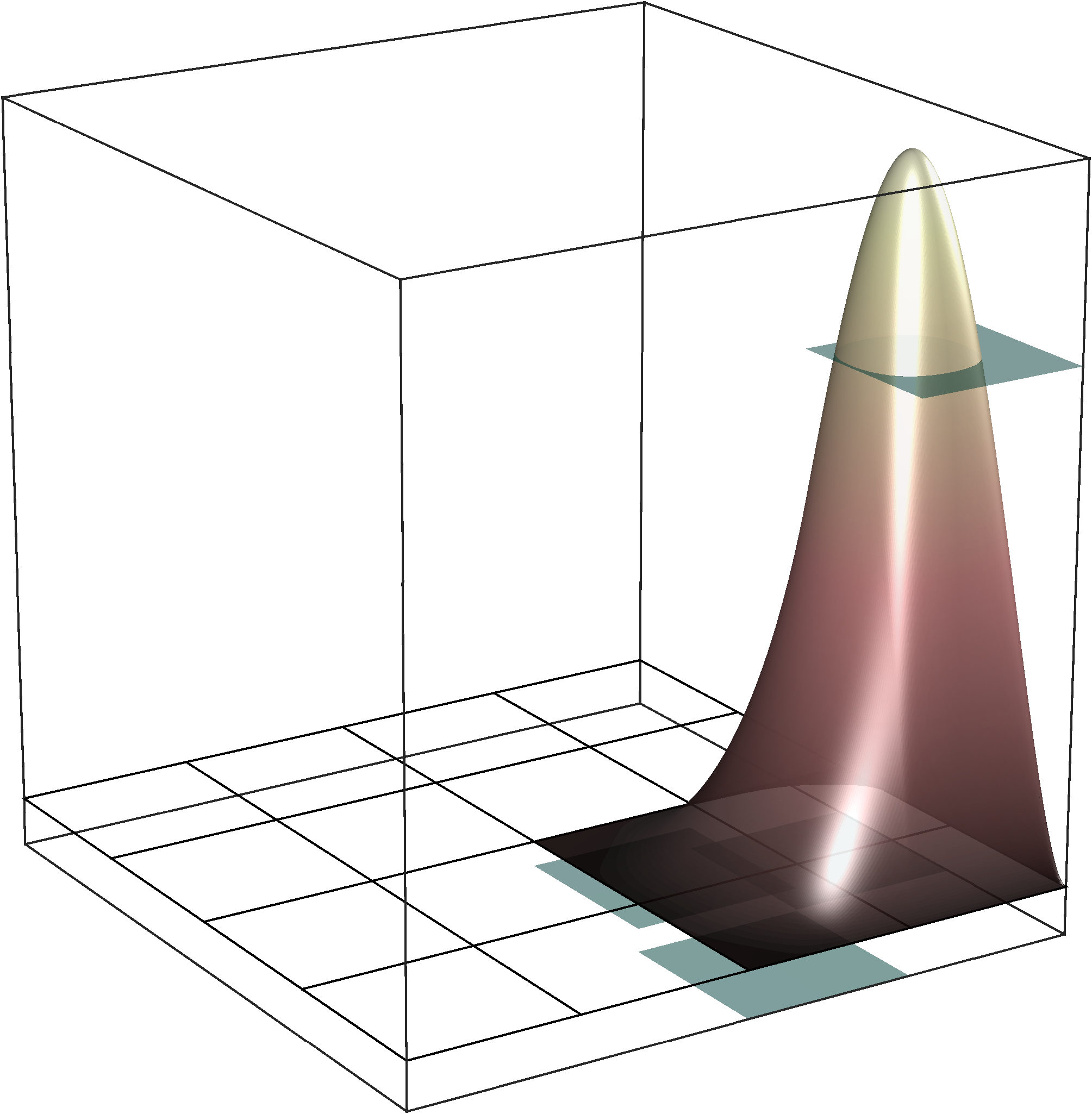}
	\hfill
	\includegraphics[width=.3\linewidth]{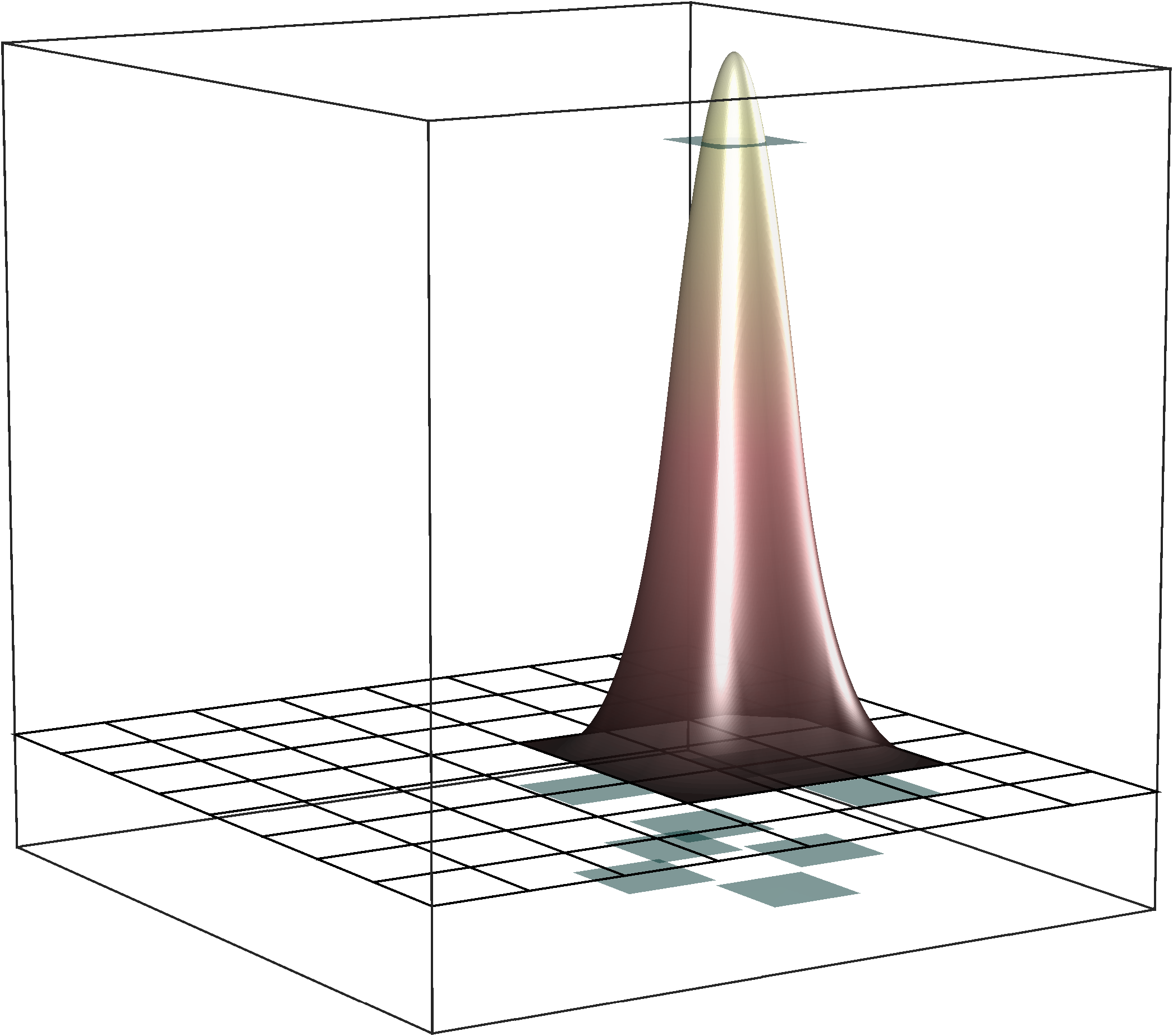}
	\hfill
	\includegraphics[width=.3\linewidth]{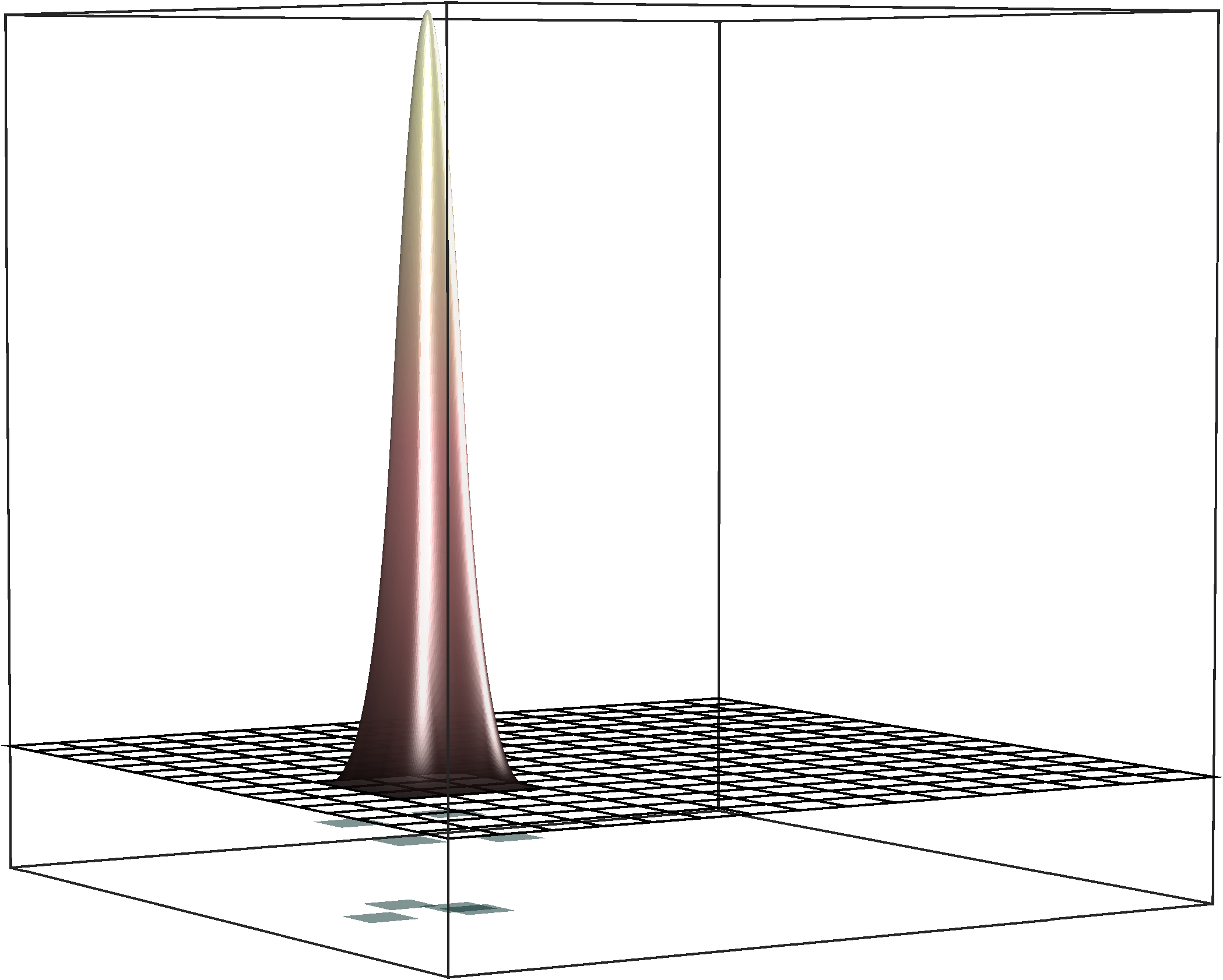}
	\caption{Super-localized basis functions obtained by the novel localization approach for $A \equiv 1$ and $\ell = 1$ in two spatial dimensions. The corresponding ($L^2$-normalized) right-hand sides $g$ are depicted in green.}
	\label{figslod2}
\end{figure}

\begin{remark}[Super-localized basis in 2d]
	In higher dimensions, the novel basis construction does not lead to a local basis. However, we will see later that the resulting localization errors decay much faster than for the their classical LOD counterparts. For an illustration of the novel basis in two dimensions, see Figure \ref{figslod2}. 
\end{remark}

\section{Super-localized numerical homogenization}\label{s:locmethod}
Using the novel localization strategy of the previous section, we turn the prototypical numerical homogenization method \eqref{e:galerkinideal} into a feasible scheme. For a fixed oversampling parameter $\ell$, define the ansatz space of the localized method as the span of the localized basis functions $\varphi^\mathrm{loc}_{\elem,\ell}$ defined in \eqref{e:patchproblem}, i.e., 
\begin{equation*}
\V_{H,\ell} \coloneqq \linh\{\varphi^\mathrm{loc}_{\elem,\ell}\,|\, \elem \in \mathcal T_H\}\subset \V.
\end{equation*}
The localized method determines the Galerkin approximation in the space $\V_{H,\ell}$, i.e., it  seeks $u_{H,\ell}\in \V_{H,\ell}$ such that, for all $v_{H,\ell}\in \V_{H,\ell}$,
\begin{equation}\label{e:wflocmethod}
a(u_{H,\ell},v_{H,\ell}) = \tspf{f}{v_{H,\ell}}{\Omega}.
\end{equation}
Let $\mathcal A_{H,\ell}^{-1}\colon  L^2(\Omega)\rightarrow \V_{H,\ell}$, $f\mapsto u_{H,\ell}$ denote the localized discrete solution operator of the above Galerkin method. 
\begin{remark}[Collocation version]\label{r:memory}
	An alternative discrete approximation, denoted also by $ u_{H,\ell}\in \V_{H,\ell}$, is given by
	\begin{equation*}
	 u_{H,\ell} = \sum_{\elem \in \mathcal T_H}c_\elem \, \varphi^\mathrm{loc}_{\elem,\ell},
	\end{equation*}
	where $(c_\elem)_{\elem \in \TH}$ are the coefficients of the expansion of $\PiH f$ in the basis functions $g_{\elem,\ell}$.
	This method is related to collocation in the sense that it enforces the PDE to hold (up to localization errors) on average in each element of $\T_H$. Similar as for the Galerkin method, this variant only requires the solution of a linear system of equations on the coarse scale. The collocation version has the advantage that the assembly of the system matrix requires only the coarse functions $g_{T,\ell}$ and does not require the computation of any inner products between the $\varphi^\mathrm{loc}_{H,\ell}$.
\end{remark}

A minimal requirement for the stability and convergence of the Galerkin method \eqref{e:wflocmethod} and its collocation variant from Remark \ref{r:memory} is that  $\{g_{\elem,\ell}\,|\,\elem \in \TH\}$ spans $\Pnull (\mathcal T_H)$ in a stable way. Numerically, this can by ensured as outlined in Appendix \ref{s:impl}. For the subsequent numerical analysis we make the following assumption.
\begin{assumption}[Riesz stability]\label{a:Rieszbasis}
The set $\{g_{\elem,\ell}\,|\,\elem \in \TH\}$ is a Riesz basis of  $\Pnull(\mathcal T_H)$, i.e., there is $C_{\mathrm{rb}}(H,\ell)>0$ depending polynomially on $H, \ell$ such that, for all $(c_\elem)_{T \in \TH}$,
	\begin{equation*}
	C^{-1}_{\mathrm{rb}}(H,\ell)\sum_{\elem \in \mathcal T_H} c_\elem^2  \leq \Big\lVert \sum_{\elem \in \mathcal T_H} c_\elem g_{\elem,\ell}\Bigl\rVert_{L^2(\Omega)}^2 \;\leq C_{\mathrm{rb}}(H,\ell) \sum_{\elem \in \mathcal T_H} c_\elem^2.
	\end{equation*}
\end{assumption}

\section{Abstract error analysis}\label{s:error}
In this section, we derive an error estimate for the localized solution operator $\mathcal A_{H,\ell}^{-1}$ introduced in Section~\ref{s:locmethod}. The error estimate is explicit in the quantity 
\begin{equation*}\label{e:qom}
\sigma = \sigma(H,\ell) \coloneqq \max_{\elem\in\TH} \sigma_\elem(H,\ell)
\end{equation*}
which reflects the maximal localization error over the coarse mesh, cf. \eqref{e:orth}.
Bounds for $\sigma$ will be presented at the end of this section and in Section~\ref{s:qo}. 

\begin{theorem}[Uniform localized operator approximation]\label{t:error}
		Let $\{g_{\elem,\ell}\,|\,\elem \in \TH\}$ be stable in the sense of Assumption \ref{a:Rieszbasis}. For any $s\in [0,1]$, the localized discrete solution operator $\mathcal A^{-1}_{H,\ell}$  approximates $\mathcal A^{-1}$ in the operator norm $\|\cdot \|_{H^s(\Omega)\rightarrow \V}$ with an algebraic rate in $H$ plus an additional localization error, i.e., there are $C,\,C^\prime>0$ independent of $H,\ell$ such that, for all $f \in H^s(\Omega)$,
	\begin{equation*}
	\begin{aligned}
	{\vnormf{\mathcal A^{-1}f - \mathcal A_{H,\ell}^{-1}f}{\Omega}}&\leq C\hspace{.62ex}\big( H\tnormf{f-\PiH f}{\Omega}  + C_{\mathrm{rb}}^{1/2}(H,\ell)\ell^{d/2}\sigma(H,\ell)\|f\|_{L^2(\Omega)}\big)\\
	&\leq  C^\prime\big(H^{1+s} + C_{\mathrm{rb}}^{1/2}(H,\ell)\ell^{d/2} \sigma(H,\ell)\big)\|f\|_{H^s(\Omega)}
	\end{aligned}
	\end{equation*}
	with $C_{\mathrm{rb}}$ from Assumption \ref{a:Rieszbasis}. 
\end{theorem}
\begin{remark}[High contrast]
	The constants $C$, $C^\prime$, $\sigma$, and $C_\mathrm{rb}$ may depend on $\alpha,\beta$ but are independent of variations of $A$. Presumably, there is a similar contrast-dependence as for the LOD. In practice, the sensitivity towards the contrast ${\beta}/{\alpha}$ is usually smaller than predicted by the theory. Subsequently, constants may depend on $\alpha,\beta$ and this dependence will not be elaborated.
\end{remark}
\begin{proof}
	Using C\'ea's Lemma (see e.g. \cite{Cia78}), we obtain that  $u_{H,\ell} = \mathcal A_{H,\ell}^{-1} f$ is the best approximation of $u = \mathcal A^{-1} f$ in $\V_{H,\ell}$ (w.r.t. the energy norm), i.e.,  for all $v_{H,\ell} \in \V_{H,\ell}$,
	\begin{align*}
	\vnormf{u-u_{H,\ell}}{\Omega}\leq \vnormf{u-v_{H,\ell}}{\Omega}.
	\end{align*}
	Setting $\tilde u = \mathcal A^{-1}\PiH f$, we obtain with the triangle inequality
		\begin{align*}
	\vnormf{u-u_{H,\ell}}{\Omega}\leq  \vnormf{u-\tilde u}{\Omega}+\vnormf{\tilde{u}-v_{H,\ell}}{\Omega}.
	\end{align*}
	For the first term, we obtain using property \eqref{e:L2proj}.
	\begin{align*}
	\vnormf{u-\tilde u}{\Omega}= \vnormf{\mathcal A^{-1}(f-\Pi_H f)}{\Omega} \leq   \|f-\Pi_Hf\|_{H^{-1}(\Omega)} 
	\leq \pi^{-1}\alpha^{-1/2}  H  \tnormf{f - \PiH f}{\Omega}.
	\end{align*}
	For the second term, we use that the prototypical method \eqref{e:galerkinideal} is exact for $\PiH f$ (Remark~\ref{r:exactness}). This enables us to represent $\tilde u$ using the ideal basis functions $\varphi_{\elem,\ell}$, i.e., 
	\begin{equation*}\label{e:u_h}
	\tilde u = \sum_{\elem \in \mathcal T_H} c_\elem\, \varphi_{\elem,\ell},
	\end{equation*}
	where $(c_\elem)_{T \in \TH}$ are the coefficients of the expansion of $\PiH f$ in the basis functions $g_{\elem,\ell}$ (Assumption \ref{a:Rieszbasis}). 
	For the particular choice \begin{equation*}
	v_{H,\ell} \coloneqq \sum_{\elem\in \mathcal T_H} c_\elem\, \varphi^\mathrm{loc}_{\elem,\ell}\in \V_{H,\ell},
	\end{equation*}
	we obtain defining $e\coloneqq\tilde u - v_{H,\ell}$ 
	\begin{align}
	\label{e:proof1}
	\vnormf{e}{\Omega}^2 = \sum_{\elem \in \mathcal T_H}c_\elem\,a( \varphi_{\elem,\ell}-\varphi^\mathrm{loc}_{\elem,\ell},e).
	\end{align}
	The definition of $\sigma$ in \eqref{e:orth} shows that
	\begin{equation}\label{e:proof2}
	\begin{aligned}
	a(\varphi_{\elem,\ell}-\varphi^\mathrm{loc}_{\elem,\ell},e) &= \tspf{g_{\elem,\ell}}{\gamma_0^{-1}\gamma_0\, e}{\mathsf N^\ell(\elem)} \leq \sigma(H,\ell)\|\gamma_0^{-1}\gamma_0\, e\|_{H^1(\mathsf N^\ell(\elem))}\\&\leq C  \sigma(H,\ell)\|e\|_{H^1(\mathsf N^\ell(\elem))}.
	\end{aligned}
	\end{equation}
	with $\gamma_{0}$ and $\gamma_{0}^{-1}$ denoting the trace and $A$-harmonic extension operators on $\mathsf N^\ell(T)$, respectively. 
	For the proof of the last inequality, we decomposed $\gamma_0^{-1}\gamma_0 e = e+e_0$ with  $e_0 \in H^1_0(\mathsf{N}^\ell(T))$ satisfying, for all $v \in H^1_0(\mathsf N^\ell(T))$, $a(e_0,v) = -a(e,v)$. The inequality then follows, using
		\begin{align*}
		\|e_0\|_{H^1(\mathsf N^\ell(T))}^2 
		 &= \tnormf{\nabla e_0}{\mathsf N^\ell(T)}^2+\tnormf{e_0}{\mathsf N^\ell(T)}^2
		\leq \alpha^{-1}\big(1+\pi^{-2}\diam (\mathsf N^\ell(T))^2\big)a(e_0,e_0)
		\\
		&\leq C\alpha^{-1} |a(e,e_0)| \leq C\beta \alpha^{-1} \|e\|_{H^1(\mathsf{N}^\ell(T))}\|e_0\|_{H^1(\mathsf{N}^\ell(T))},
	\end{align*}
where we employed Friedrichs' inequality and that the diameter of the patches $\mathsf N^\ell(T) \subset \Omega$ is bounded (recall that $\Omega$ is of unit size).

	The combination of \eqref{e:proof1}, \eqref{e:proof2}, Assumption \ref{a:Rieszbasis}, the discrete Cauchy--Schwarz inequality, the finite overlap of the patches, and \eqref{e:L2proj} yields 
	\begin{equation*}
	\begin{aligned}
	\vnormf{e}{\Omega}^2 &= \sum_{\elem \in \mathcal T_H}c_\elem\,a( \varphi_{\elem,\ell}-\varphi^\mathrm{loc}_{\elem,\ell},e)\leq  C\sigma(H,\ell)\sum_{\elem \in \mathcal T_H}c_\elem\|e\|_{H^1(\mathsf{N}^\ell(\elem))} \\& \leq C\sigma(H,\ell)\sqrt{\sum_{\elem \in \mathcal T_H}c_\elem^2}\sqrt{\sum_{\elem \in \mathcal T_H}\|e\|^2_{H^1(\mathsf{N}^\ell(\elem))}}\\&\leq C C_{\mathrm{rb}}^{1/2}(H,\ell)\ell^{d/2}\sigma(H,\ell)\vnormf{e}{\Omega} \tnormf{f}{\Omega}
	\end{aligned}
	\end{equation*}
	with $C>0$ being a generic constant independent of $H,\ell$. For the last step, we again used 
	Friedrichs' inequality. Rewriting the estimates in terms of operators, the assertion follows. 
\end{proof}

\begin{remark}[Error estimate for collocation version]
	The solution of the collocation version introduced in Remark \ref{r:memory} also satisfies the error estimate from Theorem \ref{t:error}. This is an intermediate result from the proof of Theorem \ref{t:error}, where an error estimate for the collocation version is used to bound the error of the Galerkin version.
\end{remark}

The remaining part of this section contains an existence result that is tentatively pessimistic. We show that there exists a basis that satisfies   Assumption \ref{a:Rieszbasis} and that decays exponentially. For this result, the choice of basis is equivalent to the LOD \cite{MalP20,Peterseim2021} and, hence, different from the newly proposed one. The result shows that the localization error of the novel localization strategy decays at least exponentially. 

The idea is to construct the right-hand sides $g_{\elem,\ell}$ such that applying the local inverse operator yields the respective LOD basis functions. For this purpose, we introduce non-negative bubbles  $b_\elem\in H^1_0(\elem)$ with $\Pi_H b_\elem = \one_\elem$ which are chosen such that
 \begin{align}
\label{e:est-L2-bubble-func}
\pi\|b_\elem\|_{L^2(\elem)}
\leq H\, \|\nabla b_\elem\|_{L^2(\elem)}\leq C \sqrt{ |\elem| }
\end{align}
with $C>0$ solely depending on the shape regularity of the element $T \in \TH$.
The LOD basis function which corresponds to element $\elem$ and is supported on the patch $\omega = \mathsf N^\ell(\elem)$ is then given by
\begin{equation}\label{elodbasisfun}
\varphi^\mathrm{loc} =
\varphi^{\mathrm{loc}}_{\elem,\ell}\coloneqq (1-\C^\mathrm{loc})b_\elem
\end{equation}
with $ \C^\mathrm{loc} = \C^\mathrm{loc}_{\elem,\ell}$ denoting the localized correction operator which is  defined as the $a$-orthogonal projection from $H^1_0(\omega)$ onto the closed subspace $\W= \W_{\elem,\ell} \coloneqq \{w\in H^1_0(\omega)\,|\linebreak \Pi_{H,\omega} w = 0\}$ of functions oscillating at scales unresolved by $\TH$, i.e., for all $w \in \W$,
\begin{equation*}
a(\C^\mathrm{loc} v,w) = a(v,w).
\end{equation*}
Due to the coercivity of $a$, the correction operator is  well-defined and continuous with constants independent of $H,\ell$. The LOD basis function (locally in $\omega$) possesses a $\TH$-piecewise constant right-hand side (see Lemma \ref{l:L2flux}), i.e., 
\begin{equation}\label{e:gn}
g  = g_{\elem,\ell} \coloneqq -\mathrm{div}A\nabla \varphi^\mathrm{loc}_{\elem,\ell} \in \Pnull(\T_{H,\omega})\subset L^2(\omega).
\end{equation}
The following Lemma shows that this choice yields an exponentially small value of $\sigma$ and that the $g_{\elem,\ell}$ form a stable basis of $\mathbb P^0(\TH)$ in the sense of Assumption \ref{a:Rieszbasis}. 

\begin{lemma}[Stability and exponential decay of classical LOD basis]\label{t:elod}
	Suppose that $\ell$ is chosen at least proportional to $|\log H|$. Then, there is $C>0$ independent of $H,\ell$, such that $L^2$-normalized versions of the   $g_{\elem,\ell}$  satisfy, for all $(c_\elem)_{\elem \in \TH}$, 
	\begin{equation}\label{e:rb}
	(CH^{-4})^{-1}\sum_{\elem \in \mathcal T_H}c_\elem^2  \leq  \Big\|\sum_{\elem \in \mathcal T_H}c_\elem\, g_{\elem,\ell}\Big\|_{L^2(\Omega)}^2 \leq  CH^{-4} \sum_{\elem \in \mathcal T_H}c_\elem^2.
	\end{equation}
	Furthermore, inserting the normalized $g_{\elem,\ell}$  in \eqref{e:orth}, yields 
		\begin{equation}\label{e:qo}
	\sigma(H,\ell) \leq C^\prime H^{-1}\exp(-C\ell)
	\end{equation}
	with $C^\prime >0$ independent of $H,\ell$.
\end{lemma}
\begin{proof}
	For the proof, see Appendix \ref{s:appendix}.
\end{proof}

Combining Theorem \ref{t:error} and Lemma \ref{t:elod}, one recovers an a-priori error estimate with rates as for the classical LOD based on piecewise constant finite elements; see \cite{Peterseim2021}. 

\section{Justification of super-localization}\label{s:qo}
This section theoretically justifies the super-exponential decay of the localization error by utilizing a conjecture related to the decay of Steklov eigenfunctions. This decay is in line with the numerical results  in Section \ref{s:ne}. In the remainder of this section, we consider domains in $\mathbb R^d$, $d\geq 2$, as nothing needs  to be shown in one dimension; see Remark \ref{r:1d}.

We construct approximation spaces of $Y$ defined in  \eqref{e:Y} using Steklov eigenfunctions, i.e.,  solutions to the patch eigenvalue problem 
\begin{equation*}
\mathrm{div} A\nabla \psi = 0  \text{ in }\omega,\quad
A\nabla \psi\cdot n = \lambda \psi  \text{ on }\Gamma_1 \coloneqq \partial \omega \backslash \partial \Omega,\quad 
\psi = 0  \text{ on }\Gamma_2\coloneqq\partial \omega \backslash \Gamma_1
\end{equation*}
with $n$ denoting the outer normal unit vector with respect to $\omega$.  The weak formulation of this  eigenproblem seeks eigenpairs $(\psi,\lambda) \in \V_\omega\times \mathbb R$ such that, for all $v \in \V_\omega$,
\begin{equation*}\label{e:steklov}
a_\omega(\psi,v) = \lambda \tspf{\psi}{v}{\Gamma_1}.
\end{equation*}
There exists an $L^2(\Gamma_1)$-orthonormal and $a_\omega$-orthogonal set of eigenfunctions  $\{\psi_k\,|\,k\in \mathbb{N}_0\}$  with corresponding non-negative eigenvalues $\{\lambda_k\,|\,k \in \mathbb N_0\}$ (ordered such that $0 = \lambda_0 \leq \lambda_1\leq \dots$) which is a complete subset of $Y$; see \cite{Auchmuty05}. 
Arbitrary functions in $Y$ can be expanded as 
\begin{equation}\label{e:harmproj}
u = \sum_{k = 1}^{\infty} \tspf{u}{\psi_k}{\Gamma_1}\psi_k,
\end{equation}
where the sum converges in the $H^1(\omega)$-norm. 

\begin{figure}[h]
	\includegraphics[width=.3\linewidth]{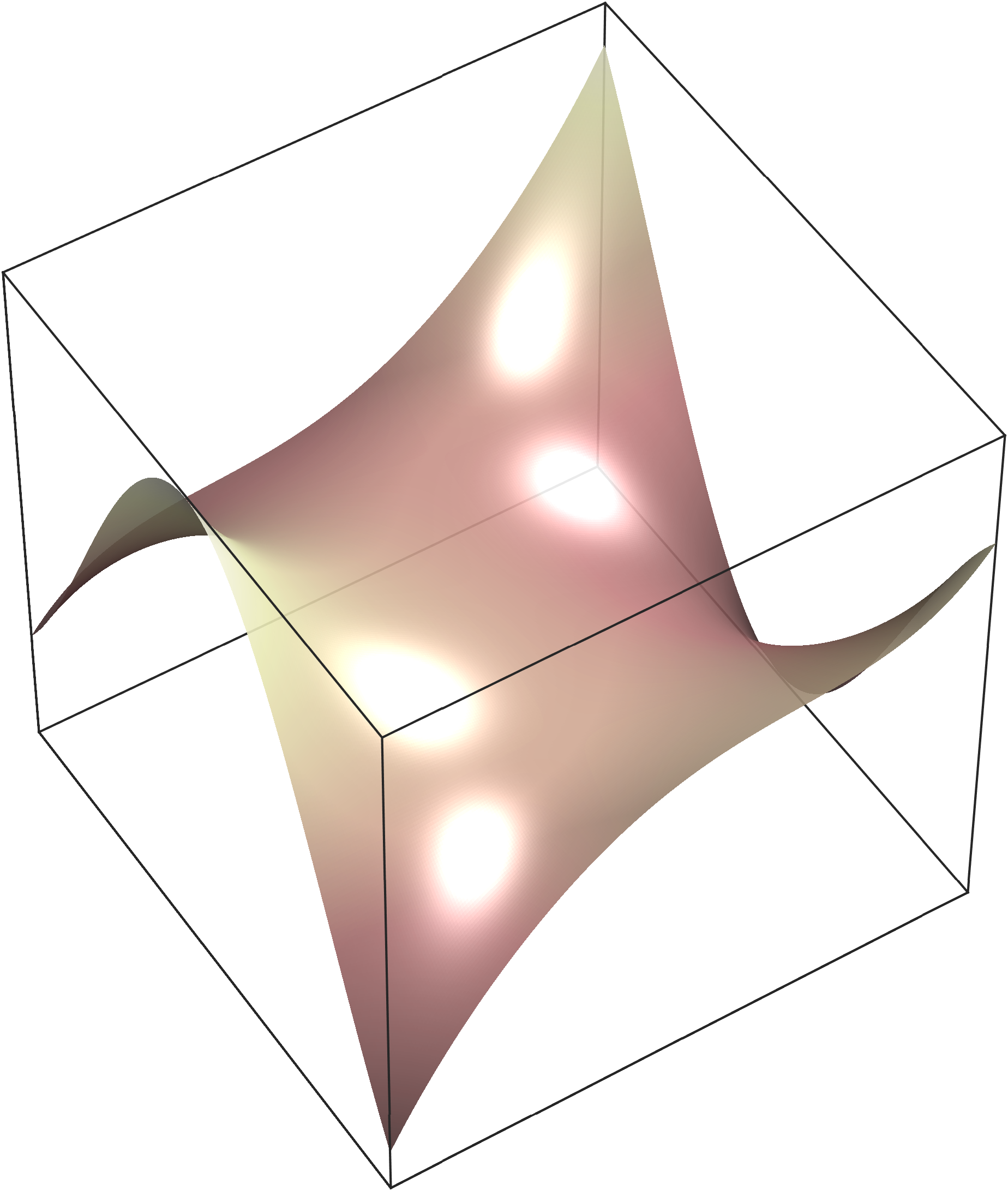}
	\includegraphics[width=.3\linewidth]{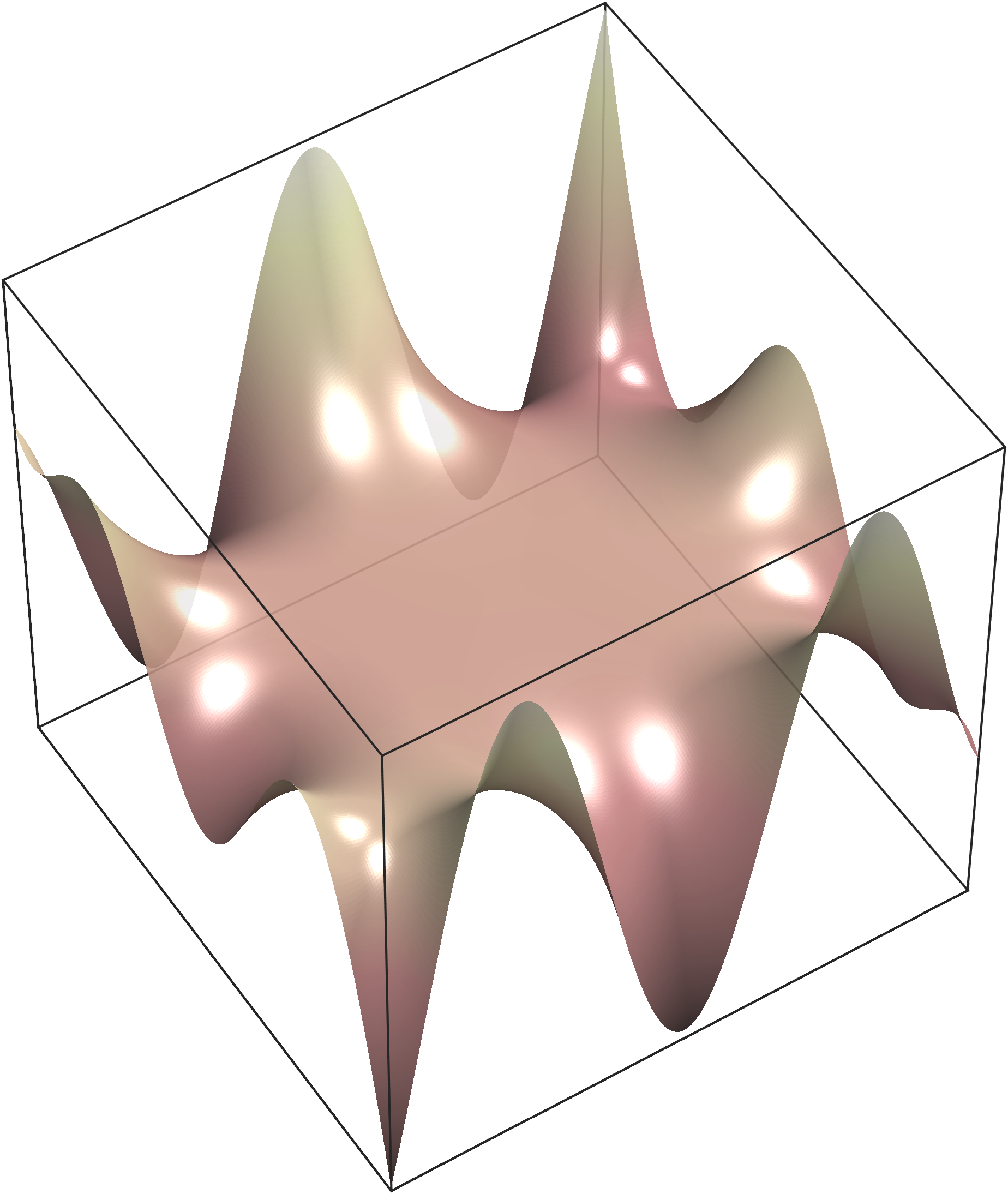}
	\includegraphics[width=.3\linewidth]{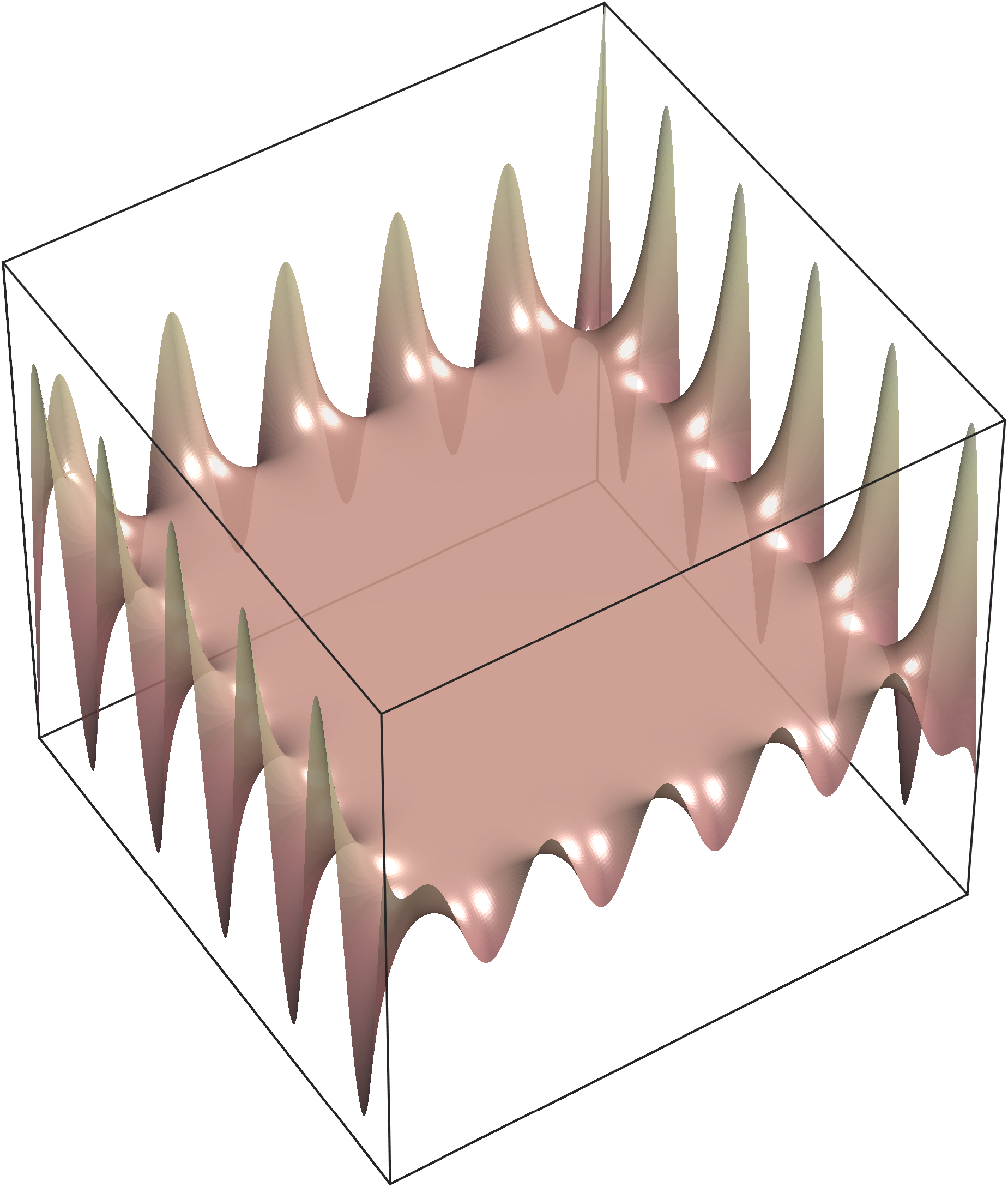}
	\caption{Illustration of Steklov eigenfunctions $\psi_6,\psi_{15},$ and $\psi_{42}$ for $\omega = (0,1)^2$, $\Gamma_1 = \partial \omega$, and $A \equiv 1$. }
	\label{figsteklov}
\end{figure}

A sufficient condition for the super-exponential decay of the localization error can be phrased in terms of certain decay properties of the Steklov eigenfunctions.

\begin{assumption}[Decay of element-averaged Steklov eigenfunctions]\label{a:stelov}
	There are $c>0$  independent of $H,\ell$ and $C_\mathrm{sd}>0$ depending polynomially on $H,\ell$ such that for all patches $\mathsf N^\ell(\elem)$ the Steklov eigenfunctions $\psi_k$ satisfy, for all $k \in \mathbb N_0$,
	\begin{equation*}
	\tnormf{\Pi_{H,\mathsf N^\ell(T)}\psi_k}{\mathsf N^{\lfloor \ell/2\rfloor}(\elem)} \leq C_{sd}(H,\ell) \exp(-c k^{\frac{1}{d-1}}).
	\end{equation*}
\end{assumption}

\begin{remark}[Decay of Steklov eigenfunctions]
	The oscillatory nature of the Steklov eigenfunctions for large $k$ presumably leads to a strong decay of the corresponding eigenfunctions; see Figure \ref{figsteklov}. For some special cases, this has been made rigorous; see e.g. \cite{Hislop2001,Polterovich15,Galkowski16}. These works require smoothness (or partially even analyticity) of the domain and the coefficients.  In \cite{Hislop2001}, a conjecture on the decay of the Steklov eigenfunctions in the interior of a domain has been formulated. This conjecture in particular implies the element-averaged version from Assumption \ref{a:stelov}.
\end{remark}

\begin{theorem}[Super-exponential localization]\label{t:epsest}
	Suppose that Assumption \ref{a:stelov} holds and let $L \in \mathbb N$ be the largest number such that, for all	$\ell \leq L$ and all patches, the number \linebreak of patch elements scales like $\ell^d$. Then, there are $c,\,C>0$ independent of $H,\ell$ such that, for all $\ell \leq L$, 
	\begin{equation*}
	\sigma(H,\ell) \leq C C_\mathrm{sd}(H,\ell)\ell^{d-\frac{d}{d-1}} H^{-1/2} \exp(-c\ell^{\frac{d}{d-1}})
	\end{equation*}
	with   $C_{\mathrm{sd}}$ from Assumption \ref{a:stelov}.
\end{theorem}
\begin{proof}
	Define $\widehat \omega \coloneqq \mathsf N^{\lfloor \ell/2\rfloor}(\elem)$ and let $N_{\widehat \omega} = \#\T_{H,{\widehat\omega}}$ denote the number of elements in $\widehat \omega$. We estimate $\sigma_\elem$ by choosing the $L^2$-normalized $g\in \Pnull(\T_{H,{\widehat\omega}})$ in \eqref{e:orth} such that it is $L^2$-orthogonal to the first $N_{\widehat\omega}-1$ Steklov eigenfunctions, i.e., 
	\begin{equation*}\label{e:choiceg}
	g \perp_{L^2(\omega)} \{\psi_1,\dots,\psi_{N_{\widehat\omega}-1}\}.
	\end{equation*} 
	If this choice is non-unique, one may require $L^2$-orthogonality to even more Steklov eigenfunctions (until $g$ is uniquely determined). This choice may not coincide with the optimal one in \eqref{e:gopt} but suffices to derive an upper bound for $\sigma_\elem$. By \eqref{e:harmproj} and Assumption \ref{a:stelov}, this choice yields, for all $v \in Y$, using 
	\begin{align*}
	\tspf{g}{v}{{\omega}} &= \Big(g\,,\,\sum_{k = N_{\widehat\omega}}^{\infty} \tspf{v}{\psi_k}{\Gamma_1}{\Pi_{H,\omega}}\psi_k\Big)_{L^2({\widehat\omega})} \leq \Big\|\sum_{k = N_{\widehat\omega} }^{\infty} \tspf{v}{\psi_k}{\Gamma_1}{\Pi_{H,\omega}}\psi_k\Big\|_{L^2({\widehat\omega})}\\
	&\leq \tnormf{v}{\Gamma_1}\sum_{k = N_{\widehat\omega}}^{\infty}\tnormf{{\Pi_{H,\omega}}\psi_k}{{\widehat\omega}}\leq C_\mathrm{sd}(H,\ell) \,\tnormf{v}{\Gamma_1}\sum_{k = N_{\widehat\omega}}^{\infty}\exp(-c k^{\frac{1}{d-1}}).
	\end{align*}
	Rewriting the last sum as a (generalized) geometric sum with the base $\theta \coloneqq \exp(-c) <1$ and estimating the sum against an integral, we get
	\begin{align*}
	\sum_{k = N_{\widehat\omega}}^{\infty}\exp(-c k^{\frac{1}{d-1}}) =  \sum_{k = N_{\widehat\omega}}^{\infty} \theta^{k^{\frac{1}{d-1}}}\leq C^\prime \int_{N_{\widehat{\omega}}}^{\infty}\theta^{x^{\frac{1}{d-1}}}\dx,
	\end{align*}
	where $C^\prime>0$ denotes a generic constant independent of $\elem,H,\ell$. 
	Using a change of variables and integrating by parts until the term $x^{d-2}$ in the product $x^{d-2}\theta^x$ has vanished, one obtains
		\begin{equation*}\label{e:x}
		\int_{N_{\widehat{\omega}}}^{\infty}\theta^{x^{\frac{1}{d-1}}}\dx = (d-1) \int_{{N_{\widehat\omega}^{\frac{1}{d-1}}}}^{\infty} x^{d-2}\theta^x \dx\leq C^\prime N_{\widehat\omega}^{\frac{d-2}{d-1}} \theta^{N_{\widehat\omega}^\frac{1}{d-1}}.
	\end{equation*}
	With the definition of $\theta$ and $N_{\widehat\omega}\approx \ell^d$, we obtain 
		\begin{align*}
	\sigma_\elem(H,\ell) &\leq \sup_{v\in Y\colon \|v\|_{H^1({\widehat\omega})} = 1}|\tspf{ g}{v}{{\widehat\omega}}|\\
	&\leq \sup_{v\in Y\colon \|v\|_{H^1({\widehat\omega})} = 1} C^\prime C_\mathrm{sd}(H,\ell)\ell^{d-\frac{d}{d-1}} \exp(-c^\prime\ell^{\frac{d}{d-1}}) \tnormf{v}{\Gamma_1}.
	\end{align*}
	introducing $c^\prime>0$ that only differs from $c$ by a constant factor. Employing the trace inequality \cite[Lemma 1.49]{dipietro12} and using the quasi-uniformity of the mesh  $\TH$, we obtain
	\begin{equation*}
	\sigma_\elem(H,\ell)\leq C^\prime C_\mathrm{sd}(H,\ell)\ell^{d-\frac{d}{d-1}} H^{-1/2} \exp(-c^\prime\ell^{\frac{d}{d-1}})
	\end{equation*}	
	Taking the maximum over all $\elem \in \TH$ yields the result. 
\end{proof}

We highlight that this proof only utilizes the decay of Steklov eigenfunctions which have an index greater or equal than the number of elements in $\mathsf N^{\lfloor{\ell}/{2}\rfloor}(T)$. Thus, provided that $\ell$ is chosen sufficiently large, a fixed number of non-decaying Steklov eigenfunctions does not pose a problem in the above proof and 
Assumption \ref{a:stelov} could be relaxed accordingly.

\begin{remark}[Theorem \ref{t:epsest} in 1d]
	Interpreting $\tfrac{d}{d-1}$ as infinity for $d=1$, Theorem \ref{t:epsest} is consistent with the earlier result that the new  basis is local in 1d (see Remark \ref{r:1d}).
\end{remark}

\section{Numerical experiments}\label{s:ne}
A practical implementation of the novel numerical homogenization method requires a fine-scale discretization for solving the patch problems  \eqref{e:patchproblem} and for the computation of the respective right-hand sides \eqref{e:gopt}. This can be done by substituting the spaces $H^1(\omega)$, $H^1_0(\omega)$, $X$, and $Y$, defined in Section \ref{s:locmethod}, by their finite element counterparts on meshes obtained by successive uniform refinement of $\T_{H,\omega}$. Given the extensive experience with the numerical analysis of fine-scale discretization for the classical LOD \cite{MalP20,Peterseim2021}, we expect that the theoretical results remain valid for this case.

The subsequent numerical experiments shall illustrate the super-localized numerical homogenization method and its theoretical properties. We consider uniform Cartesian meshes of $\Omega = (0,1)^d$,  $d \in \{2,3\}$. Henceforth, let the mesh size $H$ denote the side length of the elements instead of the diameter. 

For the numerical experiments in 2d, we choose $A$ to be a realization of the coefficient which is piecewise constant with respect to the mesh of mesh size $2^{-7}$ and which takes independent and identically distributed element values between $\alpha = 0.01$ and $\beta = 1$. 
For the numerical experiment in 3d, we use the simpler periodic coefficient 
\begin{equation*}
A(x_1,x_2,x_3) = \frac{1}{100} + \frac{99}{200}\Big(\prod_{i = 1}^{3}\sin(2^8\pi  x_i)+1\Big)
\end{equation*}
which has the same values of $\alpha $ and $\beta$. The periodicity of this coefficient on the mesh of mesh size $2^{-7}$ is exploited in the numerical implementation in order to reduce the computational cost; see Remark \ref{r:per}. Details on the implementation are provided in Appendix \ref{s:impl}.
For all numerical experiments (2d and 3d), the  local patch problems \eqref{e:patchproblem} are solved using the $\mathcal Q_1$-finite element method on Cartesian meshes of the respective patches with mesh size $h = 2^{-9}$.
The fully discrete numerical approximations (denoted by $u_{H,h,\ell}$) are obtained by the  Galerkin method \eqref{e:wflocmethod} for the experiments in 2d and by the collocation variant from Remark \ref{r:memory} for the experiments in 3d.
The energy errors are calculated with respect to a reference solution (denoted by $u_h$) which is obtained by the $\mathcal Q_1$-finite element method on the global mesh with mesh size $h = 2^{-9}$. 

\subsection*{Observation of super-exponential localization}
\begin{figure}[h]
	\includegraphics[width=.499\linewidth]{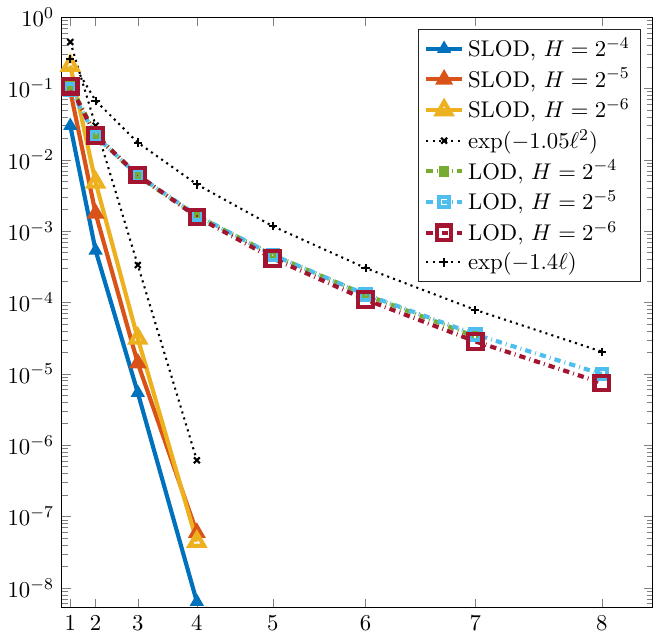}\hfill\hfill
	\includegraphics[width=.499\linewidth]{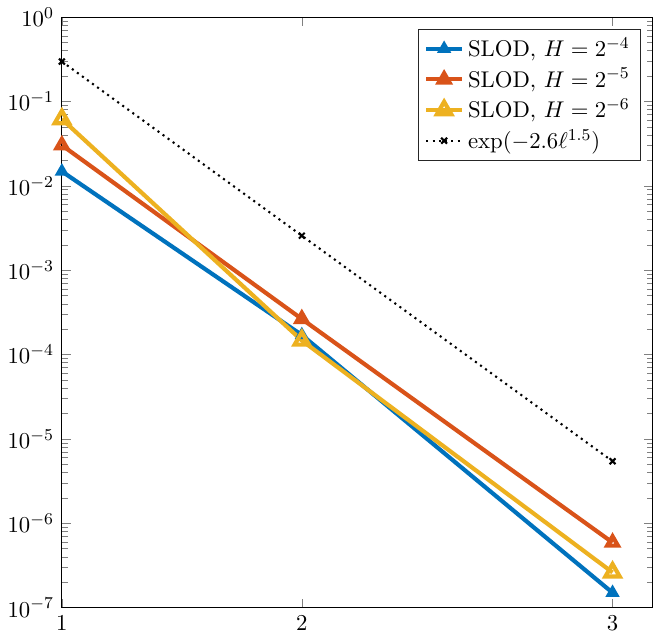}
	\caption{Plot of the energy localization errors $\vnormf{u_{H,h,\ell}-u_h}{\Omega}$ in dependence of the oversampling parameter $\ell$ of the SLOD and the LOD for the 2d experiment (left) and of the SLOD for the 3d experiment (right).}
	\label{fig1}
\end{figure}

We consider the right-hand side $f \equiv 1$, as for $f\in\Pnull(\TH)$, the error is bounded solely by the localization error (see Theorem~\ref{t:error}).  In Figure \ref{fig1}, the localization errors for the novel localization approach (referred to as SLOD) are shown  for the above (2d and 3d) numerical experiments. The localization errors are plotted for several coarse grids $\TH$ in dependence of $\ell$. For the 2d experiment, we additionally depict the localization errors of the stabilized LOD from \cite{HaPe21}. For reference, we indicate lines showing the expected rates of decay of the localization errors.

Figure \ref{fig1} numerically confirms the super-exponential decay rates of the localization errors as stated in Theorem \ref{t:epsest}. The localization error of the LOD, depicted in Figure~\ref{fig1}, decays exponentially; see e.g. \cite{MaP14,Peterseim2021}. Much higher values of $\ell$ are necessary in order to reach the accuracy level of the SLOD. 
For the 3d experiment,  the reference problem needs $2^{3\cdot 9} \approx 130$M DOF for sufficiently resolving the coefficient. By numerical homogenization, without significant loss of accuracy, it is possible to reduces the number of DOF to, e.g., $2^{3\cdot 4}\approx 4$K for $H = 2^{-4}$ which demonstrates the high potential of the novel localization technique.

\subsection*{Optimal convergence under mesh-refinement}
\begin{figure}[h]
	\includegraphics[width=.499\linewidth]{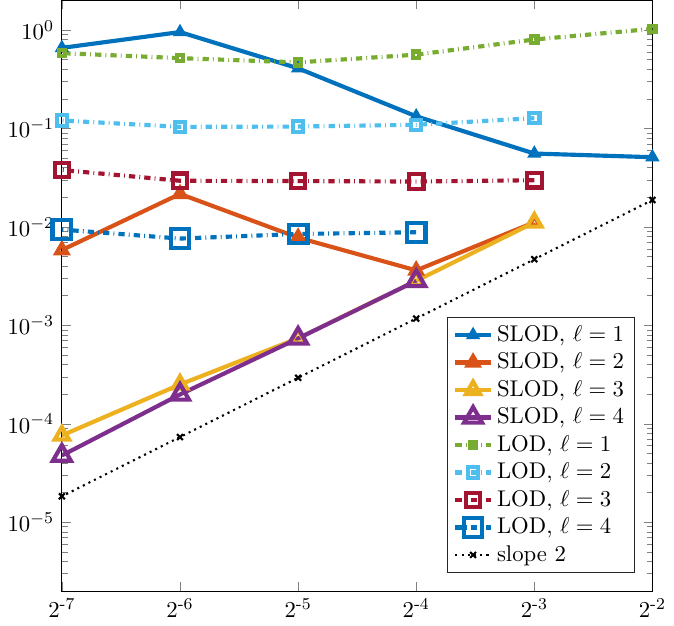}
	\caption{Plot of the energy errors $\vnormf{u_{H,h,\ell}-u_h}{\Omega}$ in dependence of the mesh size $H$ of the SLOD and the LOD for the 2d experiment.}
	\label{fig3}
\end{figure}

For demonstrating the convergence of the new method, we consider the 2d numerical example with  right-hand side
\begin{equation*}
f(x_1,x_2) = \sin(x_1)\sin(x_2).
\end{equation*}
Figure \ref{fig3} shows the energy errors for the SLOD and the LOD for different oversampling parameters as the coarse mesh $\TH$ is refined. Note that we only consider combinations of $H$, $\ell$ for which all patch-problems \eqref{e:patchproblem} are non-global.  For reference, a line of slope 2 is depicted.

\section{Conclusion}\label{s:conc}
This paper presents and analyzes a novel localization approach for elliptic multiscale problems. It constructs an approximation space with optimal algebraic rates of convergence in $H$ without pre-asymptotic effects which has  basis functions with supports of width $\mathcal O(H|\log H|^{(d-1)/d})$. A sequence of numerical experiments illustrates the significance of the novel localization strategy when compared to the so far best localization to supports of width $\mathcal O(H|\log H|)$. 

Three questions remain open and need to be addressed in future work.
The first one is the stable selection of a basis of right-hand sides $g_{\elem,\ell}$ (Assumption \ref{a:Rieszbasis}). We have presented a practical implementation that numerically solves this issue; see  Appendix \ref{s:impl}. Nevertheless, a proof or a novel implementation that yields a provably stable basis, is still missing. 
The second question concerns the super-exponential decay of the localization errors. This decay is mathematically justified in Section \ref{s:qo}, however, a rigorous proof is open. 
Lastly, in \cite{GGS12} it is conjectured that there exists a fully local basis (i.e. $\mathcal O(H^{-d})$ basis functions with supports of diameter $\mathcal O(H)$) with optimal rates of convergence without pre-asymptotic effects which is  motivated by numerical experiments in two dimensions. The present work neither proves nor disproves the existence of such a local basis, however, due to the super-exponential decay, the novel basis functions are extremely close to being local. 

Despite these open questions, the approach has certainly high potential beyond the simple model problem of this paper, as the recent works \cite{FHP21,BFP22} show.

\section*{Acknowledgments}
We would like to thank Roland Maier and Barbara Verf\"urth for inspiring discussion on the localization problem and Timo Neumeier for his valuable contributions in the early stage of the development of this paper. Furthermore, we are very grateful to Philip Freese and Christoph Zimmer for many constructive discussions and their help in improving the implementation of the method.

\appendix

\section{Proof of Lemma \ref{t:elod}}\label{s:appendix}
\begin{remark}[Tilde notation]
	Henceforth, we write $a\lesssim b$ and $c\gtrsim d$ in short for $a \leq  Cb$ and $c\geq C^\prime d$, respectively with constants $C,C^\prime>0$ independent of $H$, $\l$. 
\end{remark}

We define the $\V$-conforming $L^2$-orthogonal projection $\PitH$ onto the span of bubble functions as
\begin{equation*}
\PitH(v)\coloneqq\sum_{\elem\in\TH}\fint_\elem v \dx\, b_\elem
\end{equation*}  
with $b_\elem$ from \eqref{e:est-L2-bubble-func}. It is straightforward to prove local continuity and approximation properties on the element level, i.e., it holds, for all $T \in \TH$,
\begin{equation}\label{e:pit}
\begin{aligned}
\|\PitH  v\|_{L^2(\elem)}&\leq \|v\|_{L^2(\elem)} &&\text{ for all }v\in L^2(\elem),\\
\tnormf{\nabla \PitH v}{\elem}&\lesssim H^{-1}\tnormf{v}{\elem}\quad&& \text{ for all }v\in L^2(\elem),\\
\|v-\PitH v\|_{L^2(\elem)}&\lesssim H \|\nabla v\|_{L^2(\elem)} &&\text{ for all }v\in H^1(\elem).
\end{aligned}
\end{equation}
\begin{lemma}[Properties of LOD right-hand sides]\label{l:L2flux}
	The LOD basis function $\varphi^\mathrm{loc}$ defined in \eqref{elodbasisfun} has a right-hand side with 
	\begin{equation}\label{e:glod}
	g = g_{\elem,\ell} \coloneqq -\mathrm{div}A\nabla \varphi^\mathrm{loc} \in \Pnull(\T_{H,\omega})\subset L^2(\omega).
	\end{equation}
	Furthermore, the right-hand side decays exponentially, i.e., it holds, for all $m \in \mathbb N$, 
	\begin{equation}\label{e:decg}
	\tnormf{g}{\omega\backslash\mathsf N^m(\elem)}\lesssim H^{-1}\exp(-C m) \tnormf{g}{\omega}.
	\end{equation}  
\end{lemma}
\begin{proof}
	For proving \eqref{e:glod}, we consider the saddle-point formulation from  \cite{OwhS19,Mai20ppt} which seeks the localized LOD basis function $\varphi^{\mathrm{loc}} = \varphi_{\elem,\ell}^\mathrm{loc}$ and a corresponding Lagrange multiplier $\lambda \in \Pnull(\T_{H,\omega})$ such that
	\begin{align*}
	\begin{pmatrix}
	\mathcal A & \mathcal B^T\\
	\mathcal B & 0
	\end{pmatrix} \begin{pmatrix}
	\varphi^{\mathrm{loc}}\\\lambda
	\end{pmatrix}
	= \begin{pmatrix}
	0\\\mathbf{1}_\elem
	\end{pmatrix}
	\end{align*}
	with the patch-local operators $\mathcal A\colon  H^1_0(\omega) \rightarrow H^{-1}(\omega), u \mapsto a_\omega(u,\cdot)$ with $a_\omega$ denoting the restriction of $a$ to $\omega$, $\mathcal B\colon  H^1_0(\omega)\rightarrow \Pnull(\T_{H,\omega}), u \mapsto \Pi_{H,\omega}  u$, and its transpose $\mathcal B^T\colon  \Pnull(\T_{H,\omega})\rightarrow H^{-1}(\omega)$ defined for all $p \in \Pnull(\T_{H,\omega})$ and $v \in H^1_0(\omega)$ by 
	$\langle \mathcal B^Tp,v\rangle_{H^{-1}(\omega)\times H^1_0(\omega)} \coloneqq \tspf{p}{v}{\omega}.$
	Next, we show that $\mathcal S\colon  \Pnull(\T_{H,\omega})\rightarrow \Pnull(\T_{H,\omega}), \,p\mapsto  (\mathcal B \mathcal A^{-1}\mathcal B^T)p$ is invertible by applying the Lax--Milgram lemma. We have that  $\mathcal S$ is coercive as for arbitrary $p\in \Pnull(\TH)$ 
	\begin{align*}
	\tspf{\mathcal Sp}{p}{\omega} = \tspf{\mathcal A^{-1}\mathcal B^Tp}{\mathcal B^T p}{\omega} = a_\omega(\mathcal A^{-1}\mathcal B^T p,\mathcal A^{-1}\mathcal B^Tp)\gtrsim \|\mathcal B^T p\|_{H^{-1}(\omega)} ^2.
	\end{align*}
	The estimate can be continued by using the properties of the bubble operator \eqref{e:pit} together with the inverse estimate $\tnormf{\nabla \PitH v}{\omega} \lesssim H^{-1} \tnormf{\PitH v}{\omega}$ for bubble functions 
	\begin{align*}
	H\tnormf{p}{\omega} &= H\sup_{v\in H^1_0(\omega)}\frac{\tspf{\mathcal B^T p}{v}{\omega}}{\tnormf{v}{\omega}} \lesssim  H\sup_{v\in H^1_0(\omega)}\frac{\tspf{\mathcal B^T p}{\PitH v}{\omega}}{\tnormf{\PitH v}{\omega}}\\
	&\lesssim \sup_{v\in H^1_0(\omega)}\frac{\tspf{\mathcal B^T p}{\PitH v}{\omega}}{\tnormf{\nabla \PitH v}{\omega}} \leq  \|\mathcal B^T p\|_{H^{-1}(\omega)}.
	\end{align*}
	Joining the estimates proves the coercivity of $\mathcal S$. The continuity of $\mathcal S$ can be seen as follows
	\begin{equation*}
	\tspf{\mathcal S p}{q}{\omega} = a_\omega(\mathcal A^{-1} \mathcal B^T p,\mathcal A^{-1} \mathcal B^T q)\lesssim  \|\mathcal B^T p\|_{H^{-1}(\omega)}\|\mathcal B^T q\|_{H^{-1}(\omega)} \lesssim \tnormf{p}{\omega}\tnormf{q}{\omega}.
	\end{equation*}
	The Lax--Milgram lemma then yields that $\mathcal S^{-1}$ exists and
	\begin{equation*}
	\tnormf{\mathcal S^{-1} p}{\omega} \lesssim H^{-2} \tnormf{p}{\omega}.
	\end{equation*}

	As  $\varphi^{\mathrm{loc}}$ satisfies 
	\begin{equation*}
	\mathcal A \varphi^{\mathrm{loc}} = \mathcal B^T\mathcal S^{-1}\mathbf{1}_\elem,
	\end{equation*}
	it is an immediate consequence that $\varphi^{\mathrm{loc}}$ solves for all $v \in H^1_0(\omega)$
	\begin{equation*}
	a_\omega(\varphi^{\mathrm{loc}},v) = \tspf{\mathcal S^{-1}\mathbf{1}_\elem}{v}{\omega}
	\end{equation*}
	and thus $g \coloneqq -\mathrm{div} A \nabla \varphi^{\mathrm{loc}} = S^{-1}\mathbf{1}_\elem\in \Pnull(\T_{H,\omega})\subset L^2(\omega)$.
	
	For the proof of \eqref{e:decg}, we utilize the connection of $\mathcal S^{-1}$ to LOD theory. Using \eqref{elodbasisfun}, one obtains for $p,q\in \Pnull(\T_{H,\omega})$
	\begin{align*}\label{e:decaylod}
	\tspf{\mathcal S^{-1} p}{q}{\omega} &= \langle \mathcal A \mathcal A^{-1}\mathcal B^T\mathcal S^{-1} p,\mathcal A^{-1}\mathcal B^T\mathcal S^{-1}q\rangle_{H^{-1}(\omega)\times H^1_0(\omega)} \\
	&= a_\omega ((1-\C^\mathrm{loc})\PitH p,(1-\C^\mathrm{loc})\PitH q).
	\end{align*}	
	
	Defining $\bar g \coloneqq g \mathbf{1}_{\Omega\backslash\mathsf{N}^m(\elem)}$, we obtain for the LOD right-hand side $g$ 
	\begin{align*} 
	\tnormf{g}{\omega\backslash \mathsf{N}^m(\elem)}^2 = \tspf{\mathcal S^{-1} \,\mathbf{1}_\elem}{\bar g}{\omega} = a_\omega((1-\C^\mathrm{loc})\PitH \mathbf{1}_\elem,(1-\C^\mathrm{loc})\PitH \bar g).
	\end{align*}
Utilizing the decay result \cite[Theorem 3.15]{Peterseim2021}, the assertion
	\begin{equation*}
	\tnormf{g}{\omega\backslash \mathsf N^m(\elem)}\lesssim H^{-1}\exp(-C m)\tnormf{g}{\omega}
	\end{equation*}
	can be concluded using standard cut-off arguments from LOD theory.
\end{proof}

\begin{proof}[Proof of Lemma \ref{t:elod}]
	We begin proving assertion \eqref{e:qo}. 
Inserting 
	 $g$ from \eqref{e:gn} in \eqref{e:orth} and using the stability result for $\gamma_0^{-1}\gamma_0$ as derived in the proof of Theorem \ref{t:error}, we obtain
	\begin{equation*}
	\sigma_\elem(H,\ell) \lesssim  \frac{1}{\tnormf{g}{\omega}} \sup_{v\in \V_\omega\colon \|v\|_{H^1(\omega)} = 1}|\tspf{g}{\gamma_0^{-1}\gamma_0\,v}{\omega}|,
	\end{equation*}
	where we account for possibly non-normalized $g$ by dividing by its norm. 
	We choose a finite element cut-off function $\eta\in W^{1,\infty}(\omega)$  with $\|\nabla \eta \|_{L^\infty(\omega)}\lesssim H^{-1}$ such that 
	\begin{equation*}
	\eta  \equiv 0 \;  \text{ in } \mathsf N^{\ell-1}(\elem),\quad
	\eta \equiv 1 \;  \text{ on } \Gamma_1,\quad
	0\leq \eta \leq 1 \;  \text{ in } \omega\backslash \mathsf N^{\ell-1}(\elem).
	\end{equation*}
	Using $\gamma_{0}\,v = \gamma_{0}\,\eta v$ ($\eta \equiv 1$ on $\Gamma_1$), we can write for $v \in \V_\omega$
	\begin{align*}\label{e:est}
	-\tspf{g}{\gamma_0^{-1}\gamma_0 \,v}{\omega} = -\tspf{g}{\eta v}{\omega} + a_\omega(\varphi^\mathrm{loc},\eta v)
	\end{align*}
	with $\varphi^\mathrm{loc}$ from \eqref{elodbasisfun}. 
	Estimating the first term, yields
	\begin{align*}
	|\tspf{g}{\eta v}{\omega}| \leq \tnormf{g}{\omega\backslash \mathsf N^{\ell-1}(\elem)}\tnormf{v}{\omega} \lesssim H^{-1}\exp(-C\ell)\tnormf{g}{\omega}\|v\|_{H^1(\omega)},
	\end{align*}
	where we used Lemma \ref{l:L2flux}. For the second term, we get using the decay of the LOD basis function (c.f. \cite[Theorem 3.15]{Peterseim2021})
	\begin{align*}
	|a_\omega(\varphi^\mathrm{loc},\eta v)| &\lesssim \tnormf{\nabla \varphi^\mathrm{loc}}{\omega\backslash \mathsf N^{\ell-1}(\elem)}\big(\tnormf{v \nabla \eta}{R} + \tnormf{\eta \nabla v}{\omega\backslash\mathsf N^{\ell-1}(\elem)}\big)\\
	&\lesssim H^{-1} \exp(-C  \ell)\tnormf{\nabla \varphi^\mathrm{loc}}{\omega}\tnormf{\nabla v}{\omega} \\
	&\lesssim H^{-1} \exp(-C  \ell) \tnormf{g}{\omega}\|v\|_{H^1(\omega)}.
	\end{align*}
	Taking the maximum over all elements $\elem \in \TH$ concludes the proof.
	
	For proving the Riesz property \eqref{e:rb}, we first examine the linear map between the characteristic functions and the ideal LOD right-hand sides more closely. Similar as in the proof of Lemma \ref{l:L2flux}, the ideal LOD basis functions and its corresponding (non-normalized) global right-hand sides are determined by
	\begin{equation*}
	\varphi_\elem \coloneqq \mathcal A^{-1}\mathcal B^T g_\elem \quad \text{and}\quad g_\elem \coloneqq \mathcal S^{-1} \mathbf 1_\elem
	\end{equation*}
	with $\mathcal S \coloneqq \mathcal B \mathcal A^{-1} \mathcal B^T$, where $\mathcal A$ and  $\mathcal B$ are now the global counterparts of the operators in the proof of Lemma \ref{l:L2flux}.
	Following the proof of Lemma \ref{l:L2flux}, one can prove that $\mathcal S^{-1}$ (here $\mathcal S$ also denotes a global operator) exists and 
	\begin{equation*}
	\tnormf{\mathcal S^{-1} p}{\Omega} \lesssim H^{-2} \tnormf{p}{\Omega}.
	\end{equation*}
	
	Second, we show that that the difference between the global $g_\elem$ and its localization $g_{\elem,\ell}$ is exponentially small with respect to the $L^2$-norm. We obtain
	\begin{align*}
	\tnormf{g_\elem-g_{\elem,\ell}}{\Omega}^2 = \tnormf{g_\elem-g_{\elem,\ell}}{\omega}^2 + \tnormf{g_\elem-g_{\elem,\ell}}{\Omega\backslash\omega}^2.
	\end{align*}	
	Using \cite[Theorem 3.19]{Peterseim2021}, the first term can be estimated as follows
	\begin{align*}
	\tnormf{g_\elem-g_{\elem,\ell}}{\omega}^2 &= \tspf{g_\elem-g_{\elem,\ell}}{\PitH|_\omega (g_\elem-g_{\elem,\ell})}{\omega}\\ &= a_\omega( \varphi_{\elem}\vert_\omega-\varphi_{\elem,\ell}^\mathrm{loc},\PitH|_\omega (g_\elem-g_{\elem,\ell}))\\
	&\lesssim H^{-1}\exp(-C\ell)\tnormf{\nabla \varphi_{\elem}}{\Omega}\tnormf{g_\elem-g_{\elem,\ell}}{\omega}\\
	&\lesssim H^{-1}\exp(-C\ell)\tnormf{g_\elem}{\Omega}\tnormf{g_\elem-g_{\elem,\ell}}{\omega},
	\end{align*}
	where $\varphi_{\elem,\ell}^\mathrm{loc}$ denotes the localization of $\varphi_\elem$.
	For the second term, we use that $g_\elem$ decays exponentially away from $\elem$. This can be proven similarly as the local result from Lemma \ref{l:L2flux}. Using the two estimates, we conclude
	\begin{align}\label{e:gerr}
	\tnormf{g_\elem-g_{\elem,\ell}}{\Omega}\lesssim H^{-1}\exp(-C\ell) \tnormf{g_\elem}{\Omega}.
	\end{align}
	
	Henceforth, let  $g_{\elem,\ell}$ and  $g_\elem$ be $L^2$-normalized. Next, we derive the Riesz stability result for the global right-hand sides. Using the  continuity estimates for $\mathcal S^{-1}$ and $\mathcal S$, one can show the desired Riesz property in the global case
	\begin{equation}\label{e:Rideal}
		H^4\sum_{\elem \in \mathcal T_H}c_\elem^2 \lesssim  \Big\|\sum_{\elem \in \mathcal T_H}c_\elem g_\elem\Big\|_{L^2(\Omega)}^2 \lesssim H^{-4} \sum_{\elem \in \mathcal T_H}c_\elem^2.
	\end{equation}
	Lastly, estimate \eqref{e:gerr} can be reformulated for the $L^2$-normalized right-hand sides as 
	\begin{equation*}
	\tnormf{g_\elem-g_{\elem,\ell}}{\Omega}\lesssim H^{-1}\exp(-C\ell) .
	\end{equation*}
	Assuming that $\ell \gtrsim |\log H|$ and using the previous estimate, one can show that the  localized right-hand sides $g_{\elem,\ell}$ form a Riesz basis similarly as \eqref{e:Rideal}. 
\end{proof}
\section{Practical implementation}\label{s:impl}

We shall finally discuss some aspects of the stable and efficient implementation of the novel localization strategy introduced in Sections \ref{s:locapprox} and \ref{s:locmethod}. The implementation involves a random sampling approach for approximating the space $Y$ from \eqref{e:Y}, similar to \cite{smetana,Chen20}. Approximating the space $Y$ by harmonic polynomials \cite{BaL11} or Steklov eigenfunctions \cite{Ma21} would also be possible.

For the sake of simplicity, we only consider $\Omega \subset \mathbb R^d$, $d\geq 2$. 
An efficient implementation of the algorithm requires, on the one hand, that the right-hand sides $g_{\elem,\ell}$ can be computed with minimal communication between the patches and, on the other hand, that $\{g_{\elem,\ell}\,|\,\elem \in \TH\}$ is a stable basis of $\Pnull(\TH)$ in the sense of Assumption \ref{a:Rieszbasis}. 
Stability issues occur whenever there are nested groups of patches  close to $\partial \Omega$. 

\begin{figure}[h]
	\includegraphics[height=.22\textwidth]{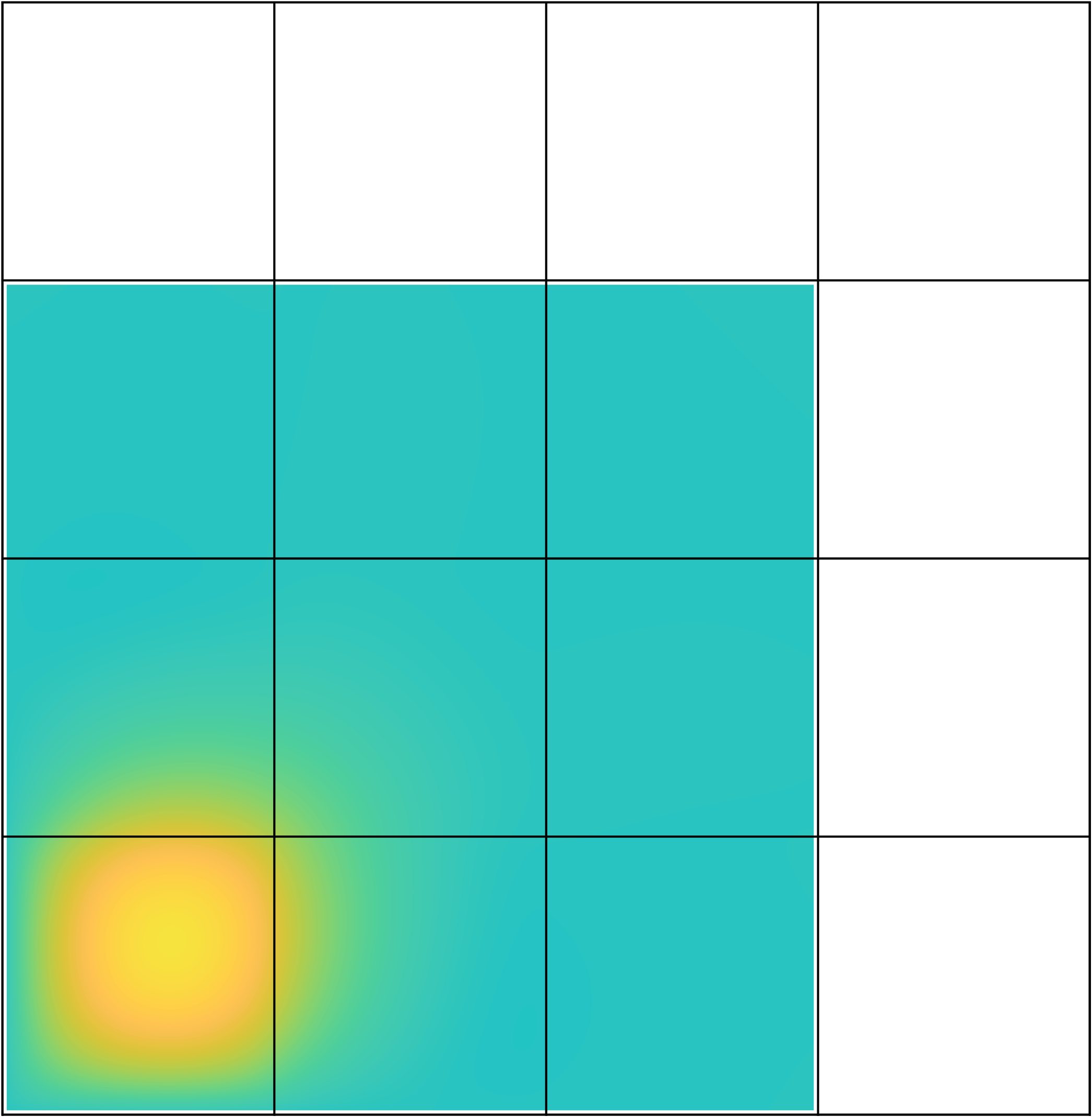}\hspace{1.4ex}\includegraphics[height=.22\textwidth]{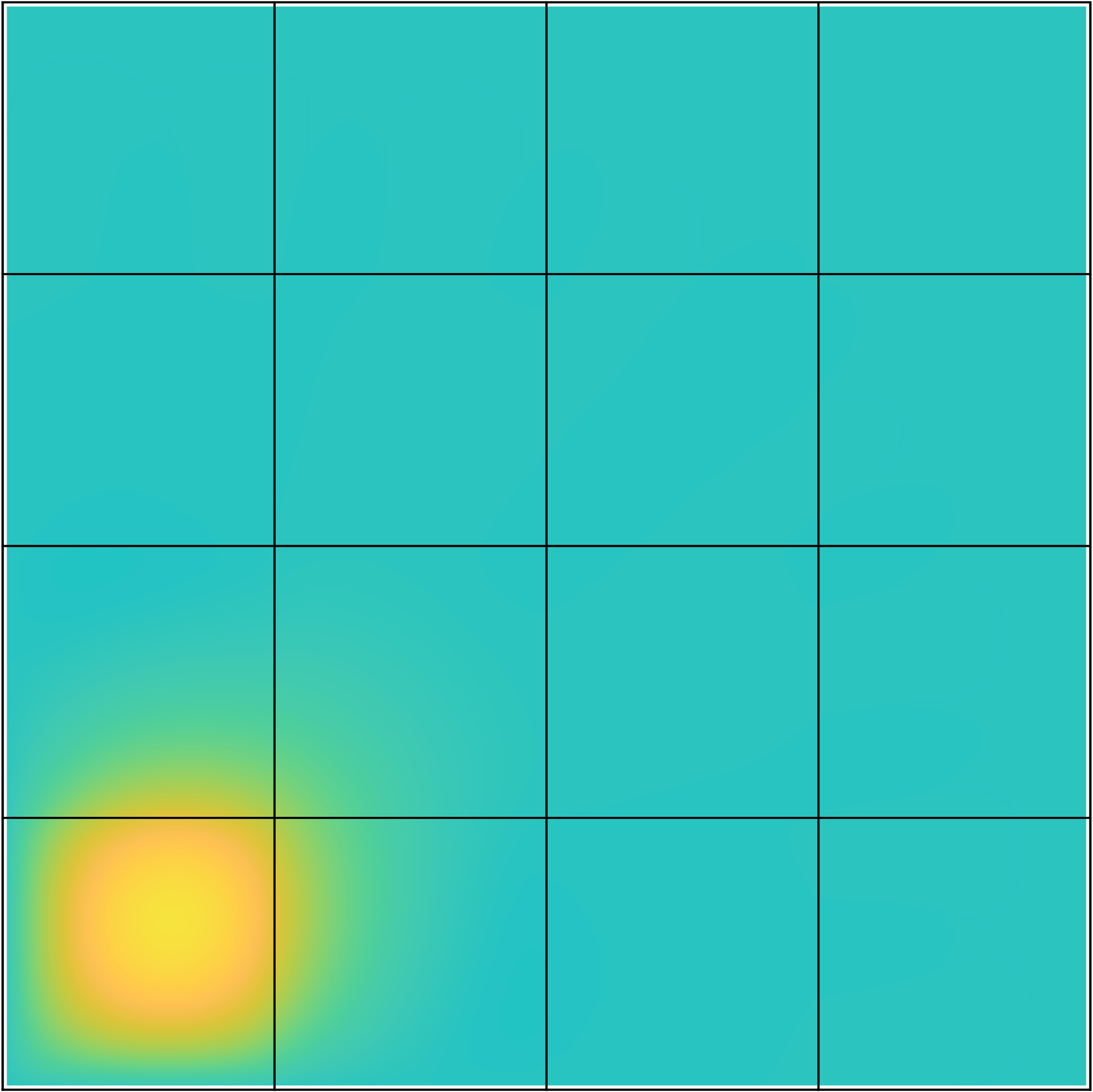}\hspace{.5ex}
	\includegraphics[height=.22\textwidth]{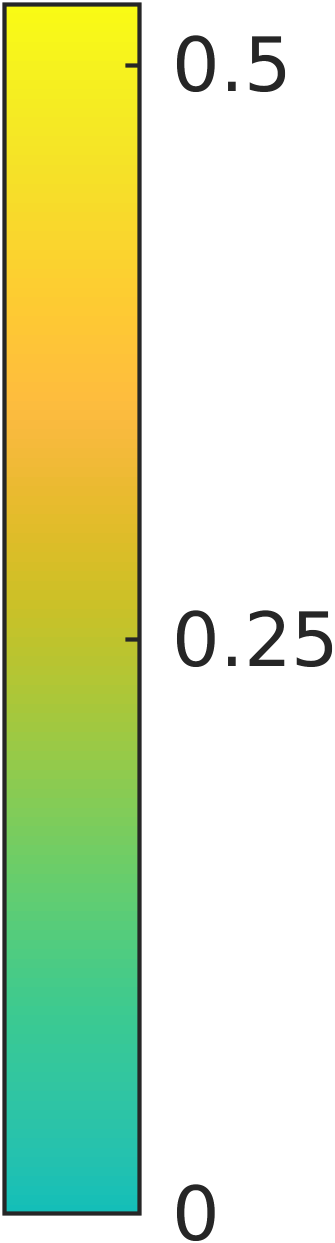}
	\caption{Functions obtained by solving the patch problems \eqref{e:patchproblem} for the right-hand sides associated to the smallest singular values on the patches corresponding to elements $(1,1)$ (left) and $(2,2)$ (right) for $A \equiv 1$.} 
	\label{fig11}
\end{figure}

For illustration purposes, consider $\Omega = (0,1)^2$ with a Cartesian mesh and $\ell = 2$. The patches corresponding to the elements with positions\footnote{Here, the position is a vector in $\{1,\dots,H^{-1}\}^2$ with first (resp. second) component  determining the location in $x$-direction (resp. $y$-direction). The numbering is so that $(1,1)$ is the bottom left element.} $(1,1),\,(1,2),$ and $(2,1)$ are contained in the patch corresponding to element $(2,2)$. This nesting leads to an unfavorable choice of basis functions, e.g., the function with right-hand side associated to the smallest singular value on the patch of element $(2,2)$ is not the basis functions one would expect for this element, but rather the one of element $(1,1)$. Indeed, it nearly coincides with the function with right-hand side associated to the smallest singular value on the patch of element $(1,1)$; see Figure \ref{fig11}. Similar problems occur for elements $(1,2)$ and $(2,1)$.

\begin{figure}[h]
	\includegraphics[width=.22\textwidth]{pics/basisill5.png}\hfill
	\includegraphics[width=.22\textwidth]{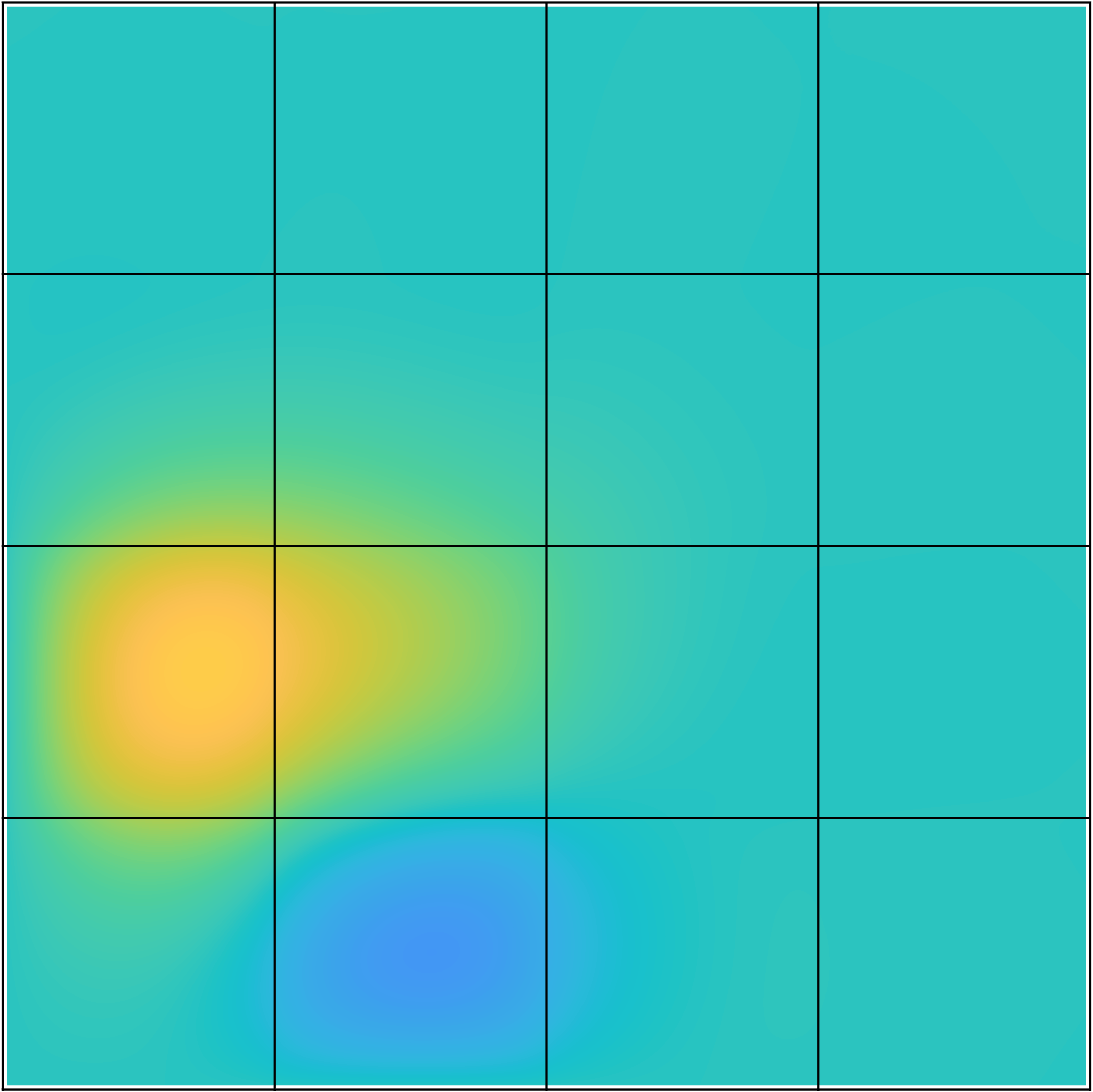}\hfill
	\includegraphics[width=.22\textwidth]{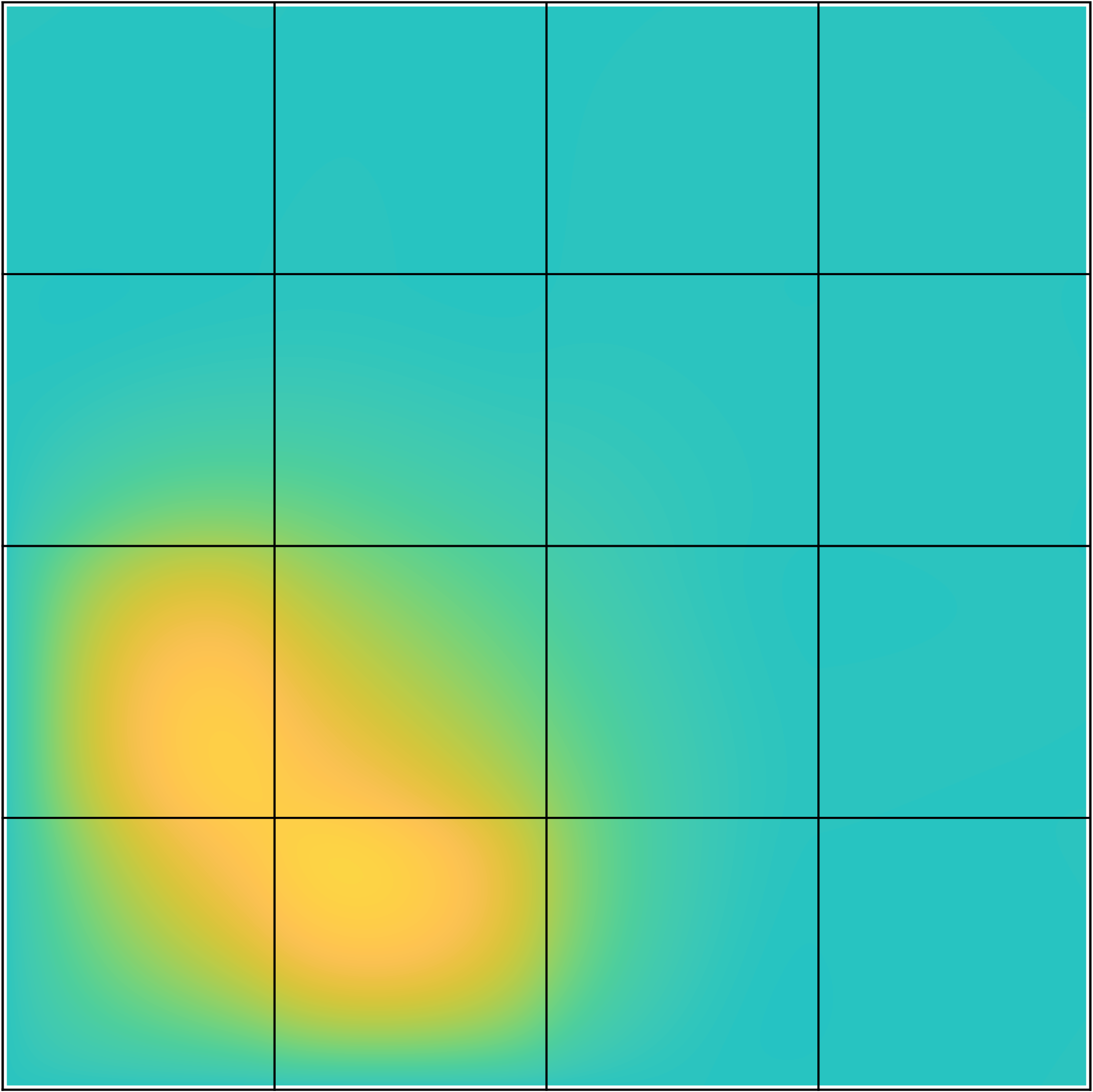}\hfill
	\includegraphics[width=.22\textwidth]{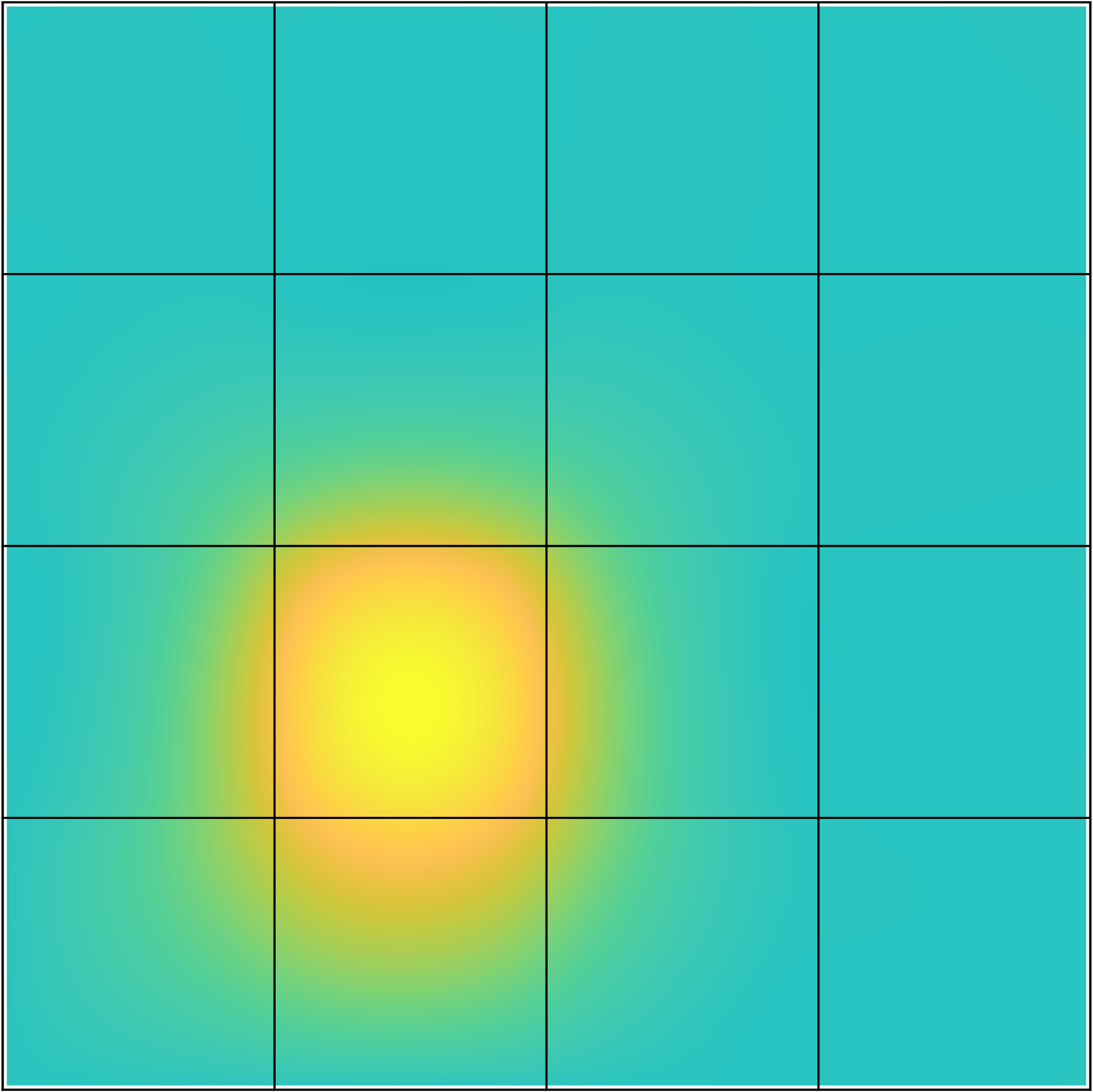}\hfill
	\includegraphics[height=.22\textwidth]{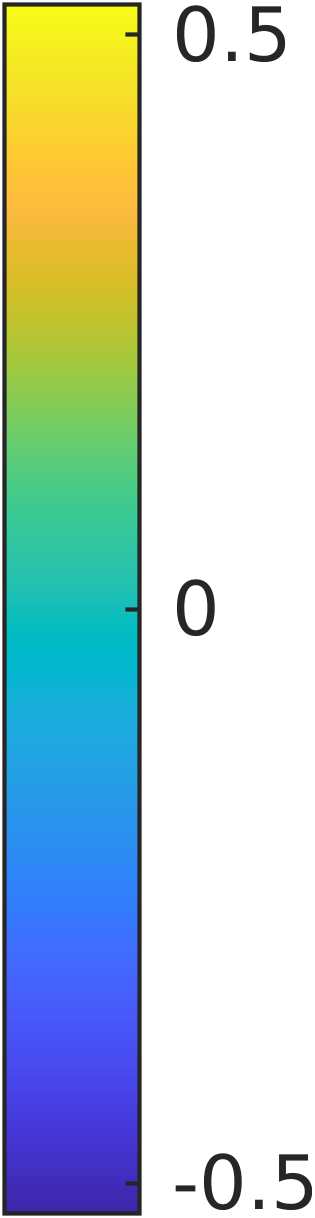}
	\caption{Functions obtained by solving the  patch problem \eqref{e:patchproblem} for the right-hand sides corresponding to the four smallest singular values (singular values increase from left to right) on the patch of element $(2,2)$ for $A \equiv 1$.}
	\label{fig22}
\end{figure}

This problem can be cured by only considering the (largest) patch of element $(2,2)$ and calculating the functions with right-hand sides associated to the four smallest singular values. For a plot of these functions, see Figure \ref{fig22}. Using this strategy, it is possible to compute four linear independent basis functions. This basis spans a space containing the functions one would intuitively expect as basis functions for the elements.

This procedure can be easily generalized. First, we identify groups of patches for which the basis functions are computed simultaneously. This can be done as follows: 
All patches $\mathsf N^\ell(\elem)$ with $\elem$ having distance $\geq \ell$ layers to $\partial \Omega$ can be considered separately. 
Patches, where $\elem$ has exactly a distance of layers $\ell-1$ to $\partial \Omega$, are considered as representatives of a corresponding group of patches.

Finally, every remaining patch $\mathsf N^\ell(T)$, i.e. $T$ has a distance of $<\ell-1$ layers to $\partial \Omega$, is uniquely assigned to a group. Here, the representative patch of this group is a superset of the patch $\mathsf N^\ell(T)$. Let the number of patches in a group be denoted by $k$. 

For all representatives $\omega = \mathsf N^\ell(\elem)$, we continue discretizing the patch in order to numerically approximate the $k$ right-hand sides. 
Denote $N \coloneqq \#\T_{H,\omega}$ and let $T_{h,\omega}$ be a fine mesh of $\omega$ obtained by uniform successive refinement of $\T_{H,\omega}$. We denote the nodes of the fine mesh $\T_{h,\omega}$ by $z_j$ with $j \in \mathcal N\coloneqq \{1,\dots,n\}$. We assume that the numbering of the nodes is such that we can disjointly decompose $\mathcal N$ into indices corresponding to nodes on $\Gamma_1\coloneqq \partial \omega \backslash \partial \Omega$, $\Gamma_2\coloneqq \partial \omega \cap \partial \Omega$, and $\omega\backslash\partial \omega$, respectively, i.e., $
\mathcal N = \mathcal N_{1} \cup \mathcal N_{2} \cup \mathcal N_0 $
with subsets $\mathcal N_{1}\coloneqq \{1,\dots,n_{1}\}$, $\mathcal N_{2} \coloneqq \{n_{1}+1,\dots.,n_{1}+n_{2}\}$, and $\mathcal N_0 \coloneqq \{n_{1}+n_{2}+1,\dots,n\}$. Furthermore, by $\{\phi_j\,|\,j \in \mathcal N\}$, we denote the hat functions associated to nodes in $\mathcal N$, satisfying $\phi_j(z_i) = \delta_{ij}$ for all $i,j \in \mathcal N$.  

We approximate the space $Y$ from \eqref{e:Y} by calculating the discrete $A$-harmonic extensions of sample boundary data on $\Gamma_1$. Denoting with $m \in \mathbb N$ the number of samples, we define $\mathbf S_{1} \in \mathbb R^{n_{1} \times m}$ to be the  matrix having sample vectors of length $n_{1}$ as columns. The entries in each column correspond to the nodal values of sample finite element boundary functions on $\Gamma_1$. 
  We use  independent and uniformly distributed nodal values within the interval $[-1,1]$.
  For the numerical experiments in Section \ref{s:ne}, the  number of samples $m$ is chosen to be $5N$.
  Let $\mathbf S^{n\times m}$ be the matrix whose columns are coordinate vectors of the discrete $A$-harmonic extensions of the  columns of $\mathbf S_{1}$. We can calculate the matrix $\mathbf S$ as 
\begin{equation*}
\mathbf S = \begin{pmatrix}
\mathbf S_{1}\\ \mathbf 0\\\mathbf S_0
\end{pmatrix}\quad \text{ with }\quad  
\mathbf S_0 \coloneqq -\mathbf A_0^{-1}\mathbf A_{1}\mathbf S_{1},
\end{equation*}
$\mathbf A_0 \coloneqq (a_\omega(\phi_j,\phi_i))_{i,j \in \mathcal N_0}$, $\mathbf A_{1} \coloneqq (a_\omega(\phi_j,\phi_i))_{i \in \mathcal N_{1}, j \in \mathcal N_0}$ and $\mathbf 0$ denoting the $n_{2}\times m$ zero matrix.
Next, the SVD of the matrix $\mathbf X \coloneqq \mathbf P (\mathbf S^T \mathbf K \mathbf S)^{-1/2} \in \mathbb R^{N\times m}$ is computed, where $\mathbf P\in \mathbb R^{N\times m}$ is the column-wise application of the $L^2$-projection onto the characteristic functions $\{\one_\elem\,|\,\elem\in \T_{H,\omega}\}$  to $\mathbf S$. The matrix $\mathbf K \in \mathbb R^{n\times n}$ is the sum of the stiffness matrix and the mass matrix with respect to  $\T_{h,\omega}$. 
It shall be noted that the term $(\mathbf S^T \mathbf K \mathbf S)^{-1/2}$ guarantees that the right singular vectors of $\mathbf X$ represent a set of $H^1(\omega)$-orthonormal functions. Nevertheless, in practice, it seems reasonable to apply the SVD directly to $\mathbf P$, i.e.,  $(\mathbf S^T \mathbf K \mathbf S)^{-1/2}$ does not need to be computed. In all our numerical experiments, we could hardly see a difference with regard to the choice of basis and the resulting errors.

The reduced SVD of $\mathbf X$ reads 
\begin{equation*}
\mathbf X = \mathbf G \mathbf \Sigma \mathbf H^T
\end{equation*}
with $\mathbf G \in \mathbb R^{N\times N}$, $\mathbf \Sigma \in \mathbb R^{N\times m}$, and $\mathbf H \in \mathbb R^{m\times m}$. The last $k$ columns of $\mathbf G$ (assuming a non-increasing ordering of the singular values) are henceforth denoted by $\mathbf g_1,\dots,\mathbf g_k \in \mathbb R^N$. The $\mathbf g_j$ are coordinate vectors (w.r.t. the basis $\{\one_\elem\,|\,\elem\in \TH\}$) of  right-hand sides $g_{\elem,\ell}$ with $\elem$ in the group of patches corresponding to representative $\omega$.

Having computed the right-hand sides $g_{\elem,\ell}$, it is not difficult to compute the fully discretized counterpart of $\varphi^\mathrm{loc}_{\elem,\ell}$ as the solution of a discretized version of \eqref{e:patchproblem}. 

\begin{remark}[Exploiting periodic coefficients]\label{r:per}
	If the coefficient $A$ is periodic with respect to $\TH$, the computational costs can be reduced significantly.  In this case, the above procedure only needs to be performed for $\mathcal O(\ell^d)$ reference patches. The basis functions can then be obtained by translation; see \cite{GaP15}.
\end{remark}

\begin{remark}[Complexity of basis computation]
Suppose that, for some $s \in \mathbb N$, we approximate the space $Y$ by $m = s\,\# \T_{H,\omega} \in \mathcal O(s\ell^d)$ samples of discrete $A$-harmonic functions on the mesh $\T_{h,\omega}$. 
Each sample requires the solution of a boundary value problem on the patch mesh with $\mathcal O((\tfrac{\ell H}{h})^d)$ degrees of freedom. 
Next, the matrix product $\mathbf S^T \mathbf K \mathbf S$ needs to be calculated which has a complexity of $\mathcal O((\tfrac{\ell H}{h})^d(s\ell^d)^2)$. Furthermore, a SVD and a matrix square root and inverse of matrices with sizes of order $\ell^d\times s\ell^d$  and $s\ell^d\times s\ell^d$, respectively, need to be computed. Finally, another boundary value problem on the patch mesh needs to be solved in order to compute the actual basis function.
\end{remark}

\bibliographystyle{alpha}
\bibliography{bib}
\end{document}